\documentclass[12pt]{amsart}
\usepackage{amssymb, tensor, fullpage, textcomp, pb-diagram, wasysym, wrapfig, comment, xspace, leftindex, xcolor, appendix, tensor, relsize}
\usepackage[Symbol]{upgreek}
\usepackage[mathscr]{euscript}
\usepackage[shortlabels]{enumitem}
\usepackage{tikz-cd}
\usepackage[all]{xy}

\let\oldtocsection=\tocsection

\let\oldtocsubsection=\tocsubsection

\let\oldtocsubsubsection=\tocsubsubsection

\renewcommand{\tocsection}[2]{\hspace{0em}\oldtocsection{#1}{#2}}
\renewcommand{\tocsubsection}[2]{\hspace{1em}\oldtocsubsection{#1}{#2}}
\renewcommand{\tocsubsubsection}[2]{\hspace{2em}\oldtocsubsubsection{#1}{#2}}

\DeclareFontFamily{U}{fsy}{}
\DeclareFontShape{U}{fsy}{m}{n}{<->s*[.9]psyr}{}
\DeclareSymbolFont{der@m}{U}{fsy}{m}{n}
\DeclareMathSymbol{\der}{\mathord}{der@m}{182}

\newcommand{\on}{\mathbf{On}}
\newcommand{\no}{\mathbf{No}}

\newcommand{\mex}{\mathrm{mex}}
\newcommand{\CB}{\mathrm{CR}}

\newcommand{\pac}{\mathrm{PAC}}
\newcommand{\prc}{\mathrm{PRC}}

\newcommand{\dlo}{\mathrm{DLO}}

\newcommand{\dprk}{\operatorname{dp}}
\newcommand{\opd}{\operatorname{opd}}

\DeclareSymbolFont{imag@m}{OT1}{cmr}{m}{ui}
\DeclareMathSymbol{\imag}{\mathord}{imag@m}{105}

\newtheorem{theorem}{Theorem}[section]
\newtheorem*{theorem*}{Theorem}

\newtheorem*{qst*}{Question}

\newtheorem{proposition}[theorem]{Proposition}
\newtheorem*{proposition*}{Proposition}
\newtheorem*{fact*}{Fact}

\newtheorem*{Claim*}{Claim}
\newtheorem{fact}[theorem]{Fact}

\newtheorem{lemma}[theorem]{Lemma}
\newtheorem{corollary}[theorem]{Corollary}
\newtheorem{conj}[theorem]{Conjecture}

\newtheorem*{thmA}{Theorem A}

\newtheorem*{thmB}{Theorem B}
\newtheorem*{thmC}{Theorem C}
\newtheorem*{thmD}{Theorem D}

\theoremstyle{definition}

\theoremstyle{remark}

\newcommand{\monster}{\boldsymbol{\Sa M}}

\newcommand{\monsterset}{\boldsymbol{M}}

\newcommand{\tp}{\operatorname{tp}}
\newcommand{\qftp}{\operatorname{qftp}}

\newcommand{\Val}{\mathrm{Val}}

\newcommand{\rgoup}{(\R;+,<)}

\newcommand{\rfield}{(\R;+,\times)}

\newcommand{\Th}{\mathrm{Th}}
\ProvideTextCommandDefault{\cprime}{(U+042C)}

\newcommand{\Fraisse}{Fra\"iss\'e\xspace}
\newcommand{\Erdos}{Erd\H{o}s\xspace}

\newcommand{\chara}{\operatorname{Char}}

\newcommand{\nip}{\mathrm{NIP}}
\newcommand{\ip}{\mathrm{IP}}

\newcommand{\Cal}[1]{\ensuremath{\mathcal{#1}}}
\newcommand{\Sa}[1]{\ensuremath{\mathscr{#1}}}

\newcommand{\age}[1]{\ensuremath{\textup{Age{#1}}}}

\newcommand{\app}{\approx_\lambda}
\newcommand{\ru}{\mathrm{RU}}
\newcommand{\ruo}{\mathrm{RU}_{\Sa O}}
\newcommand{\rum}{\mathrm{RU}_{\Sa M}}
\newcommand{\mr}{\mathrm{RM}}
\newcommand{\cra}{\mathrm{CR}}
\newcommand{\mrm}{\mathrm{RM}_{\Sa M}}

\newcommand{\Z}{\mathbb{Z}}
\newcommand{\N}{\mathbb{N}}

\newcommand{\Q}{\mathbb{Q}}
\newcommand{\R}{\mathbb{R}}

\begin{document}
\title[]{Trace definability I:\\ preservation and characterizations}

\author{Erik Walsberg}
\email{erik.walsberg@gmail.com}
\maketitle
\begin{abstract}
We introduce a notion of weak definability of first order structures, show that various classification-theoretic properties are or are not preserved under it, and that the properties which are preserved can also be characterized in terms of it.
\end{abstract}

\section*{Introduction}
A theory is unstable if and only if some formula can encode arbitrary linear orders in its models.
Thus, while an unstable structure need not define an infinite linear order outright, it does ``weakly define'' one in a precise sense.
This raises a question: what does ``weak definability" mean in general?
The aim of this paper is to develop a general theory of this form of weak definability and to show that it is a useful organizing principle in model theory.

\medskip
Our formalization is the notion of {\bf trace definability}. 
Let $\Sa M$ and $\Sa N$ be first-order structures.
If $A \subseteq N^n$, a subset of $A$ is called a {\bf trace} of an $\Sa N$-definable set if it is of the form
$Y \cap A$ for some $\Sa N$-definable set $Y \subseteq N^n$. We say that
$M$ is {\bf trace definable} in $\Sa N$ if, up to isomorphism, $M$ is a
subset of some $N^n$, and every $\Sa M$-definable subset of every $M^m$ is a trace of an $\Sa N$-definable subset of $N^{mn}$.
If we can take $n = 1$ then we say that $\Sa M$ {\bf trace embeds} into $\Sa N$.
In Section~\ref{section:ts} we show that trace embeddability is equivalent to embeddability between language free structures in a precise sense.
Trace definability is obtained from the usual definition of definability of an isomorphic copy by dropping the requirement that the domain of the copy be definable.
It provides a flexible framework for detecting structures that are not  literally definable, but whose definable sets are nevertheless visible as traces of definable sets in the ambient structure.
We show that trace definability is robust under natural
constructions and that several dividing lines can be
reformulated in terms of it.

\medskip
We will see that interpretability implies trace definability, so trace definabilty is better thought of as a weak notion of interpretability.
There are natural examples of structures which trace define but do not interpret certain kinds of algebraic structures, e.g. infinite groups or fields, see the examples at the end of this introduction.
Some of these examples will be treated in the sequel to this paper.

\medskip
Trace definability extends to theories in an obvious way.
Let $T, T^*$ be complete consistent first order theories.
Then $T$ trace defines $\Sa M$ when some model of $T$ trace defines $\Sa M$ and $T$ trace defines $T^*$ when $T$ trace defines every (equivalently: some) model of $T^*$.
At this point the first claim of the following should be more or less clear.

\begin{thmA}
\hspace{.000000000000000000000000000000000000000000000001cm}
\begin{enumerate}[leftmargin=*]
\item $T$ is unstable if and only if $T$ trace defines $(\Q; <)$ if and only if $T$ trace defines some infinite linear order.
\item $T$ is $\ip$ if and only if $T$ trace defines the \Erdos-Rado graph if and only if $T$ trace defines the generic binary relation.
\end{enumerate}
\end{thmA}

We say that two theories are trace equivalent if each trace defines the other and two structures are trace equivalent when their theories are.
So $\Sa M$ and $\Sa N$ are trace equivalent if the classes of theories that trace define $\Sa M, 
\Sa N$ agree.
Hence the \Erdos-Rado graph and the generic binary relation are trace equivalent.
The generic binary relation interprets the \Erdos-Rado graph but the \Erdos-Rado graph does not interpret the generic binary relation (we will not prove the latter here).
Matching our intuition, the \Erdos-Rado graph, modulo trace equivalence, has a canonical place as the least complicated $\ip$ structure.

\medskip
Note that Theorem~A(1) also shows that many classification-theoretic properties such as simplicity and $\mathrm{NSOP}$ are {\it not} preserved under trace definability.
In fact these properties have no trace-theoretic consequences whatsoever.
A theory is {\bf trace maximal} if it trace defines every structure and a structure is trace maximal when its theory is.
In general, if $P$ is a property of structures which is preserved under interpretations and elementary equivalences then either $P$ is preserved under trace definability or there is a trace maximal structure with $P$.
For many $P$, this follows from Theorem~B.
Recall that pseudofinite fields are supersimple and hence $\mathrm{NSOP}$.

\begin{thmB}
\hspace{.0000001cm}
\begin{enumerate}[leftmargin=*]
\item Any pseudofinite field is trace maximal.
More generally any pseudo real closed field which is not real or separably closed is trace maximal.
\item A theory $T$ is trace maximal if and only if there is $\Sa M \models T$ and an infinite subset $A$ of some $M^m$ such that every subset of every $A^n$ is a trace of an $\Sa M$-definable set.
\end{enumerate}
\end{thmB}

Thus trace maximality is the strongest possible form of the independence property.
It could also reasonably be called ``$\omega$-$\ip$".
Trace definability is therefore an inherently ``nippy" subject.
Theorem~A and transitivity of trace definability together imply that stability and $\nip$ are preserved under trace definability.
(This is also immediate from the usual definitions.)
These are the first examples of an important phenomenon: properties that are preserved under trace definability can often be characterized in terms of trace definability.
Trace equivalence is a bit like homotopy equivalence, a great deal of structure is annihilated, but what remains is heightened.

\medskip
Theorem~C summarizes our other results on preservation/characterization.
Given $k \ge 1$, we say that a structure is $k$-ary if every formula is equivalent to a boolean combination of $k$-ary formulas and formulas in the language of equality.
Hence a structure is unary if and only if it is interdefinable with a structure in a unary relational language.
Unary structures may seem too simple to be interesting, but from our point of view they are rich indeed.

\begin{thmC}
Let $T$ be an arbitrary theory.
\begin{enumerate}[leftmargin=*]
\item $T$ is $\ip$ if and only if $T$ trace defines every unary structure.
\item $T$ is not totally transcendental if and only if $T$ trace defines every unary structure in a countable language if and only if $T$ trace defines the unary relational structure $\Sa C$ with domain the Cantor set and relations defining every clopen subset.
\item $T$ is unstable if and only if $T$ trace defines every unary relational structure whose unary relations form a chain under containment if and only if $T$ trace defines every unary relational structure whose unary relations form a tree under containment.
\item $T$ is not superstable if and only if $T$ trace defines every unary relational structure whose unary relations form a well-founded tree of height $\le \omega$ under containment.
\item If $\lambda$ is a countable ordinal then $T$ has Morley rank $\ge \lambda$ if and only if the unary relational structure with domain $\omega^\lambda$ and relations defining $[0, \eta]$ for all $\eta < \omega^\lambda$ trace embeds into some $\Sa M \models T$.
Hence $T$ trace defines this structure if and only if $M^m$ has Morley rank $\ge \lambda$ for some $\Sa M \models T$ and $m \ge 1$.
\item Let $\kappa, \eta \ge 1$ be cardinals and let $\Sa B^\eta_\kappa$ be the structure with domain the set of functions $\kappa \to \eta$ and unary relations defining all sets of the form $\{ f : f(i) = j\}$ for $i < \kappa, j < \eta$.
Then $T$ has dp-rank $\ge \kappa$ if and only if $\Sa B^\eta_\kappa$ trace embeds into a model of $T$ for every $\eta$.
\item $T$ is not strongly dependent if and only if $T$ trace defines $\Sa B^\eta_\omega$ for every $\eta$.
\item Let $\Sa I_1, \ldots, \Sa I_n$ range over linear orders and let $\Sa I_\times$ be the unary relational structure with domain $I_1 \times \cdots \times I_n$ and relations defining all sets of the form $J_1 \times \cdots \times J_n$ for intervals $J_1, \ldots, J_n$.
Then $T$ has op-dimension $\ge n$ if and only if $\Sa I_\times$ trace embeds into a model of $T$ for every $\Sa I_1, \ldots, \Sa I_n$ if and only if the generic $n$-order trace embeds into a model of $T$.
\item Given $k \ge 1$, $T$ is $k$-$\ip$ if and only if $T$ trace defines the generic $(k + 1)$-ary hypergraph if and only if $T$ trace defines the generic $(k + 1)$-ary relation if and only if $T$ trace defines every $k$-ary structure.
\end{enumerate}
\end{thmC}

Hence if $\Sa M$ trace embeds into $\Sa N$ then the Morley rank, dp-rank, and op-dimension of $\Sa M$ is bounded above by the Morley rank, dp-rank, and op-dimension of $\Sa N$, respectively.
We also prove the analogous result for U-rank.
It follows that finiteness of Morley rank, dp-rank, op-dimension, and U-rank are all preserved under trace definability.


\medskip
Theorem~C shows that several stability/$\nip$ theoretic properties are ``unary" from our point of view and that $\nip$ is the weakest ``unary" property.
There is now interest in defining ``higher arity" classification-theoretic properties.
We claim that such properties should be defined in terms of trace definability of higher arity structures.
The final item of Theorem~C is a step in this direction.
In Section~\ref{section:collapse} we show that if $\Sa M$ admits quantifier elimination in a finite relational language and $\age(\Sa M)$ is a Ramsey class then a theory $T$ trace defines $\Sa M$ if and only if $T$ admits an uncollapsed indiscernible picture of $\Sa M$.

\medskip
It follows from Theorem~C(9) that the generic $(k + 1)$-ary relation is trace equivalent to the generic $(k + 1)$-hypergraph for each $k \ge 1$.
In fact, it follows from Proposition~\ref{prop:k nip} below that all $k$-$\ip$ theories admitting quantifier elimination in a finite $(k + 1)$-ary relational language are trace equivalent.
In general we expect that natural classes of structures can be classified modulo trace equivalence; as in Theorems~A and C, members of the class modulo trace equivalence will correspond to model-theoretic properties.
At present we classify countable unary theories; they correspond to stability-theoretic properties.

\begin{thmD}
Suppose that $\Sa M$ is a unary structure in a countable language.
Then exactly one of the following holds.
\begin{enumerate}[leftmargin=*]
\item $\Sa M$ is trace equivalent to the structure $\Sa C$ described in Theorem~C(2).
In this case a theory $T$ trace defines $\Sa M$ if and only if $T$ is not totally transcendental.
\item There is a unique countable ordinal $\lambda$ such that $\Sa M$ is trace equivalent to the unary relational structure with domain $\omega^{\omega^\lambda}$ and relations defining all intervals, and this $\lambda$ is the leading exponent of the Cantor normal form of the Morley rank of $\Sa M$.
In this case a theory $T$ trace defines $\Sa M$ if and only if some definable set in some model of $T$ has Morley rank $\ge \lambda$.
\end{enumerate}
Trace equivalence classes of countable unary theories form a linear order of order type $\omega_1 + 1$ under trace definability.
\end{thmD}

Along these lines, we conjecture that any infinite structure $\Sa M$ admitting quantifier elimination in a finite binary relational language is trace equivalent to either an infinite set with equality, $(\Q; <)$, or the \Erdos-Rado graph, and hence the class of theories that do not trace define $\Sa M$ is either the class of theories of finite structures, stable theories, or $\nip$ theories.

\medskip
Let $L$ be a relational language containing $\aleph_0$ relations of each arity.
Then any theory in a countable language is definitionally equivalent to an $L$-theory.
The set $\mathrm{Mod}(L)$ of $L$-structures with domain $\N$ has a well-known Polish topology.
Let $\equiv_\mathrm{tr}$ be the equivalence relation on $\mathrm{Mod}(L)$ given by trace equivalence.
It is easy to see that $\equiv_\mathrm{tr}$ is analytic.
Let $X$ be the set of $\Sa M \in \mathrm{Mod}(L)$ such that every non-unary relation in $L$ defines the empty set.
Then $X$ is a closed subset of $\mathrm{Mod}(L)$ and any unary structure in a countable language is bidefinable with some element of $X$.
By Theorem~D(2) $|X/\!\equiv_\mathrm{tr}\!|$ has cardinality $\aleph_1$.
It follows by Silver's dichotomy that $\equiv_\mathrm{tr}$ is not Borel.
Most likely $\equiv_\mathrm{tr}$ is complete analytic, but this seems difficult to prove.

\medskip
Finally, we list a couple motivational examples of non-trivial trace equivalences between ``algebraic" structures.
These example will be treated in future work; the interested reader can also take most of them as exercises.
\begin{itemize}[leftmargin=*]
\item Any $\nip$ structure is trace equivalent to its Shelah expansion.
In particular if $\Sa M$ is a linearly ordered $\nip$ structure then any expansion of $\Sa M$ by convex subsets of $M$ is trace equivalent to $\Sa M$.
This follows immediately from Shelah's theorem on externally definable sets.
This was the original motivating example.
\item $\rgoup$ is trace equivalent to $(\Z; +, <)$, in fact all non-singular ordered abelian groups are trace equivalent.
(Recall that a torsion-free abelian group $A$ is non-singular if $pA$ has finite index in $A$ for all primes $p$.)
\item $(\Z; +)$ is trace equivalent to the disjoint union of $(\Q; +)$ with $\Sa C$.
Equivalently $\Th(\Z; +)$ is the minimal theory with respect to trace definability above $\Th(\Q; +)$ which is not totally transcendental.
A torsion-free abelian group is trace equivalent to $(\Z; +)$ if and only if it is non-divisible and non-singular, equivalently if and only if it is superstable and not totally transcendental.
\item Any finite extension of $\Q_p$ is trace equivalent to $\Q_p$.
The induced structure on the set of $p$-adic balls in $\Q_p$ is trace equivalent to $\Q_p$ but does not interpret an infinite field.
\item The structure induced on roots of unity by $\rfield$ is trace equivalent to $\rfield$ but does not interpret an infinite field.
\item Bertalan Bodor has recently shown that any $\aleph_0$-stable $\aleph_0$-categorical disintegrated theory is trace equivalent to the trivial theory of an infinite set with equality.
In particular any stable theory admitting quantifier elimination in a finite relational language is trace equivalent to the trivial theory.
\end{itemize}

\subsection*{Acknowledgments}
James Hanson first showed that a theory is totally transcendental if and only if it does not trace define the unary relational structure with domain the Cantor set and unary relations defining every clopen subset~\cite{hanson-mo}.
Thanks to Chieu-Minh Tran for some useful comments.
This research was funded in part by the Austrian Science Fund (FWF) 10.55776/PAT1673125.

\section*{Conventions and background}
Throughout $m,n,k$ range over natural numbers.
All languages, structures, and theories are first order and all theories are consistent, complete, and deductively closed unless stated otherwise.
These assumptions ensure that any $L$-theory $T$ has cardinality $|L| + \aleph_0$, so if $\kappa$ is an infinite cardinal then $|T| \le \kappa$ if and only if $|L| \le \kappa$.
Unless mentioned otherwise, and ``definable" means ``first order definable, possibly with parameters".
This convention is also applied to interpretations; our interpretations are allowed to use parameters.
``Zero-definable" means ``definable without parameters".
If $\Sa M$ and $\Sa M^*$ are structures with the same domain then $\Sa M$ is a reduct of $\Sa M^*$ if every $\Sa M$-definable set is $\Sa M^*$-definable and $\Sa M$ and $\Sa M^*$ are interdefinable if each is a reduct of the other.
We let $\age(\Sa M)$ be the age of a structure $\Sa M$ in a relational language $L$, i.e. the class of finite $L$-structures that embed into $\Sa M$.
We let $\tp_{\Sa M}(a|A)$ be the type of a tuple $a$ from a structure $\Sa M$ over  $A \subseteq M$ and let $\tp_{\Sa M}(a) = \tp_{\Sa M}(a|\varnothing)$.
We let $S_n(\Sa M,A)$ be the Stone space of $n$-types over $\Sa M$ with parameters from $A \subseteq M$ and let $S_n(\Sa M) = S_n(\Sa M, \varnothing)$.
We also let $\Cal B[\Sa M]$ be the boolean algebra of zero-definable subsets of $M$.
Given an $L$-structure $\Sa M$ we let $\Sa M\!\upharpoonright\! L^*$ be the reduct of $\Sa M$ to a sublangauge $L^*$ of $L$.
We use similar notation for reducts of theories and types to sublanguages.
Given a set $A$ of parameters from some $L$-structure we let $L(A)$ be the expansion of $L$ by a constant symbol for each element of $A$.
We let $|x|$ be the length of a finite tuple $x$ of variables.
Given a structure $\Sa M$ and formula $\varphi(x,y)$ we let $\varphi(M^{|x|}, \beta)$ be the subset of $M^{|x|}$ defined by $\varphi(x,\beta)$ for any $\beta \in M^{|y|}$.

\medskip
Let $\uptau$ be a map $M \to N^n$.
The {\bf pullback} of $Y \subseteq N^{mn}$ by $\uptau$ is the set of $(a_1, \ldots, a_m) \in M^m$ such that $(\uptau(a_1), \ldots, \uptau(a_m)) \in Y$.
We often abuse notation and denote the map $M^m \to N^{mn}$ given by $(a_1,\ldots,a_m) \mapsto (\uptau(a_1),\ldots,\uptau(a_m))$ by $\uptau$ as well.

\subsection{Background on finitely homogeneous structures}
A homogeneous structure is a countably infinite structure in a relational language such that any isomorphism between finite substructures extends to an automorphism.
Fact~\ref{fact:homo} is standard, see \cite[2.1.3, 3.1.6]{macpherson-survey}.

\begin{fact}
\label{fact:homo}
Suppose that $L$ is a finite relational language and $\Sa M$ is a countable $L$-structure.
Then the following are equivalent:
\begin{enumerate}[leftmargin=*]
\item $\Sa M$ is homogeneous.
\item $\Sa M$ admits quantifier elimination.
\item $\age(\Sa M)$ is a \Fraisse class with limit $\Sa M$.
\end{enumerate}
\end{fact}

A \textbf{finitely homogeneous} structure is a homogeneous structure in a finite language.

\begin{fact}
\label{fact:isse}
Suppose that $S$ is a (necessarily incomplete) universal theory in a finite relational language and suppose that the finite models of $S$ form a \Fraisse class with limit $\Sa M$.
Then $\Th(\Sa M)$ is the model companion of $S$.
\end{fact}

\begin{proof}
By Fact~\ref{fact:homo} $\Th(\Sa M)$ has quantifier elimination.
By uniqueness of model companions it suffices to show that any model of $S$ embeds into an elementary extension of $\Sa M$.
By compactness  it suffices to show that any finite substructure of a model of $S$ embeds into $\Sa M$.
As $S$ is universal it suffices to show that every finite model of $S$ embeds into $\Sa M$.
This holds by the definition of the \Fraisse limit.
\end{proof}

If $S$ and $\Sa M$ are as in Fact~\ref{fact:isse} then call $\Sa M$ the \textit{generic model of $S$}.
For example we refer to the \Fraisse limit of the class of finite $k$-hypergraphs as the \textbf{generic $k$-hypergraph}.
Here and below a \textbf{$k$-hypergraph} $(V;E)$ is a set $V$ equipped with a symmetric $k$-ary relation $E$ such that $E(v_1,\ldots,v_k)$ implies that  $v_1, \ldots v_k$ are distinct.
We also refer to the generic graph as the \textbf{\Erdos-Rado graph}.

\subsection{Ordinals and Cantor rank}
\label{section:ordinals}
Let $\on$ be the class of ordinal numbers.
Recall that every ordinal $\eta$ may be uniquely expressed as a sum $\sum_{\lambda\in\on}\omega^\lambda n_\lambda$ where each $n_\lambda$ is a natural number and  $n_\lambda=0$ for all but finitely many $\lambda$.
(Here the sum and product are the usual sum and product of ordinals.)
This is the {\bf Cantor normal form} of $\eta$.
The {\bf natural sum} $\oplus$ of ordinals is given by adding Cantor normal forms coefficientwise.
Note that $\lambda + n = \lambda \oplus n$ for any $n$.
Given a set $A$ of ordinals we let $\mex(A)$ be  the least ordinal not in $A$.
Following \cite[Thm~4.5 and pg.~13]{gonshor} we can inductively define $\oplus$ as follows:
\begin{enumerate}
\item $\xi\oplus 0=\xi=0\oplus\xi$,
\item $\xi\oplus\zeta=\mex\{\xi\oplus\zeta^*,\xi^*\oplus\zeta : \zeta^*<\zeta, \xi^*<\xi\}$ when $\xi,\zeta>0$.
\end{enumerate}

We let $\xi\cdot\zeta$ be the natural product of ordinals $\xi,\zeta$.
This can be defined by multiplying Cantor normal forms as one would multiply polynomials, or by an inductive definition somewhat more complicated than that of $\oplus$.
We note that $n\cdot\xi$ agrees with the $n$-fold sum $\xi\oplus\cdots\oplus\xi$.

\medskip
Let $X$ be a compact Hausdorff topological space.
We recall the Cantor rank $\CB(p)$ of $p\in X$.
If $p$ is isolated then $\CB(p)=0$, if $p$ is a limit point of the set of points of rank $\ge \xi$ then $\CB(p)>\xi$, $\CB(p)$ is $\mex\{\xi : \xi<\CB(p)\}$ if this set is bounded above, and otherwise we declare $\CB(p)=\infty$.
If $\CB(p)<\infty$ then there is an open neighborhood $U$ of $p$ such that $$\CB(p)=\mex\{ \CB(p^*):p^*\in U, p^*\ne p\}.$$
Recall that $\CB(X)$ is the maximum of $\{\CB(p):p\in X\}$.
Finally, recall that $\CB(X) < \infty$ if and only if $X$ is scattered.
We declare $\infty \oplus \xi = \infty = \xi \oplus \infty = \infty \oplus \infty$ for any ordinal $\xi$.

\begin{lemma}
\label{lem:db}
Let $X$ and $Y$ be compact Hausdorff topological spaces.
Then we have $$\CB(X\times Y)=\CB(X)\oplus\CB(Y).$$
Hence  $\CB(X^n)=n\cdot\CB(X)$ for all $n$.
\end{lemma}

\begin{proof}
It suffices to prove the first claim.
If $\CB(X)=\infty$ then $X$ is not scattered, hence $X\times Y$ is not scattered, hence $\CB(X\times Y)=\infty$.
So we may suppose $\CB(X),\CB(Y)<\infty$.
It is enough to suppose that $(\alpha,\beta)\in X\times Y$ and show that $\CB(\alpha,\beta)=\CB(\alpha)\oplus\CB(\beta)$.
We apply induction on $\CB(\alpha,\beta)$.
If $\CB(\alpha,\beta)=0$ then $(\alpha,\beta)$ is isolated, hence $\alpha$ and $\beta$ are both isolated, hence $\CB(\alpha,\beta)=0=\CB(\alpha)\oplus\CB(\beta)$.
Suppose $\CB(\alpha,\beta)>0$.
By induction if $(\alpha^*,\beta^*)\in X\times Y$ and $\CB(\alpha^*,\beta^*)<\CB(\alpha,\beta)$, then $\CB(\alpha^*,\beta^*)=\CB(\alpha^*)\oplus\CB(\beta^*)$.
Now fix open neighborhoods $U\subseteq X$, $V\subseteq Y$ of $\alpha$, $\beta$, respectively, such that 
\begin{align*}
\CB(\alpha) &= \mex\{ \CB(\alpha^*) : \alpha^*\in U,\alpha^*\ne\alpha\}\\
\CB(\beta) &= \mex\{\CB(\beta^*):\beta^*\in U,\beta^*\ne\beta\}\\
\CB(\alpha,\beta) &= \mex \{\CB(\gamma) : \gamma \in U \times V, \gamma \ne (\alpha, \beta) \}.
\end{align*}
Let $\alpha^*$ range over $U$ and $\beta^*$ range over $V$.
Then we have:
\begin{align*}
\CB(\alpha,\beta) &= \mex\{ \CB(\alpha^*,\beta^*) : (\alpha^*,\beta^*)\ne (\alpha,\beta) \} \\
&= \mex\{\CB(\alpha^*)\oplus\CB(\beta^*) : (\alpha^*,\beta^*)\ne(\alpha,\beta) \} \\
&= \mex\{\CB(\alpha)\oplus\CB(\beta^*),\CB(\alpha^*)\oplus\CB(\beta): \alpha^*\ne\alpha,\beta^*\ne\beta\} \\
&=\CB(\alpha)\oplus\CB(\beta).
\end{align*}
Here the first equality holds by choice of $U$ and $V$, the second holds by induction, the third holds by choice of $U,V$ and monotonicity properties of $\oplus$, and the fourth holds by choice of $U,V$ and the definition of $\oplus$.
\end{proof}

If $\CB(X)<\infty$ then there are only finitely many $p\in X$ with $\CB(p)=\CB(X)$, the number of such $p$ is the \textit{Cantor degree} of $X$.
Fact~\ref{fact:mazsier} is a theorem Mazurkiewicz-Sierpi\'nski~\cite{Mazurkiewicz1920}.

\begin{fact}
\label{fact:mazsier}
Any countable compact Hausdorff space of Cantor rank $\lambda$ and Cantor degree $d$ is homeomorphic to the order topology on $\omega^\lambda d + 1 = (d \cdot \omega^\lambda) + 1$.
\end{fact}

Let $\mathrm{CD}(X)$ be the Cantor degree of $X$.

\begin{fact}
\label{fact:db}
Suppose that $X$ and $Y$ are second countable Stone spaces and we either have $\CB(X)>\CB(Y)$ or $\CB(X) = \CB(Y) < \infty$ and $\mathrm{CD}(X)\ge\mathrm{CD}(Y)$.
Then there is a continuous embedding $Y \hookrightarrow X$.
Hence $\omega^\lambda + 1$ continuously embeds into $X$ if $\cra(X) \ge \lambda$.
\end{fact}

The case when $\CB(X) < \infty$ follows easily from Fact~\ref{fact:mazsier}.
If $\CB(X) = \infty$ then the Cantor set continuously embeds into $X$ and we apply the fact that any second countable Stone space continuously embeds into the Cantor set.

\section{Basic facts about trace definability}\label{section:def}

A \textbf{trace embedding} $\uptau\colon \Sa M \hookrightarrow \Sa N$ is a map $\uptau\colon M \to N$ such that every $\Sa M$-definable subset of every $M^m$ is a pullback of an $\Sa N$-definable subset of $N^m$ by $\uptau$.
Consideration of the graph of equality shows that any trace embedding is injective.
It follows that $\Sa M$ is trace embeddable in $\Sa N$ if and only if we have the following up to isomorphism: $M \subseteq N$ and every $\Sa M$-definable subset of every $M^m$ is a trace of an $\Sa N$-definable subset of $N^m$.
Trace embeddings are clearly closed under composition and hence form a category.

\begin{lemma}\label{lem:b1}
Any elementary embedding is a trace embedding.
\end{lemma}

Lemma~\ref{lem:b1} is immediate from the definition.
Proposition~\ref{prop:qe-trace} follows from the observation that if $\Sa M$ is a substructure of $\Sa N$ then any quantifier-free definable subset of $M^m$ is a trace of a quantifier-free definable subset of $N^m$.

\begin{proposition}
\label{prop:qe-trace}
Suppose that $\Sa M$ is an $L$-structure which admits quantifier elimination and $\uptau\colon\Sa M \hookrightarrow \Sa N$ is an embedding of $L$-structures.
Then $\uptau$ is a trace embedding.
\end{proposition}


Given a structure $\Sa M$ and $n \ge 1$, we let $\Sa M[n]$ be the structure induced by $\Sa M$ on $M^n$, i.e. $\Sa M[n]$ has domain $M^n$ and a subset of $(M^{n})^m = M^{mn}$ is $\Sa M[n]$-definable if and only if it is $\Sa M$-definable.
(So $\Sa M[n]$ is well-defined up to interdefinability.)
A \textbf{trace definition} $\Sa M \rightsquigarrow \Sa N$ is a trace embedding $\Sa M \hookrightarrow \Sa N[n]$ for some $n \ge 1$.
More directly, a trace definition $\Sa M\rightsquigarrow \Sa N$ is a map $\uptau\colon M \to N^n$ for some $n \ge 1$ such that every $\Sa M$-definable subset of every $M^m$ is a pullback of an $\Sa N$-definable subset of $N^{mn}$ by $\uptau$.
We say that $\Sa M$ is trace definable in $\Sa N$ \textit{via} $\uptau$.
We also call $n$ the arity of $\uptau$.
Note that a trace definition must be injective.
Hence $\Sa M$ is trace definable in $\Sa N$ if and only if we have the following up to isomorphism: $M \subseteq N^n$ for some $n$ and every $\Sa M$-definable subset of every $M^m$ is a trace of an $\Sa N$-definable subset of $N^{mn}$.

\medskip
Suppose that $\uptau_i \colon \Sa M_i \rightsquigarrow \Sa M_{i + 1}$ is an $m_i$-ary trace definition for $i = 1,2$.
Let $\uptau_i = (\uptau^1_i,\ldots,\uptau^{m_i}_i)$ for $i = 1,2$, so each $\uptau^j_i$ is a function $M_i \to M_{i + 1}$.
We define the composition $\uptau_2\uptau_1$ to be the trace definition $\Sa M_1 \rightsquigarrow \Sa M_3$ given by
$$(\uptau_2\uptau_1)(a) = (\uptau^1_2(\uptau^1_1(a)),\ldots,\uptau^{m_2}_2(\uptau^1_1(a)),\ldots,\uptau^{1}_2(\uptau^{m_1}_1(a)),\ldots,\uptau^{m_2}_2(\uptau^{m_1}_1(a))).$$
It is easy to see that $\uptau_2\uptau_1$ is indeed a trace definition and that this notion of composition defines a category of trace definitions.

\medskip
Let $T,T^*$ be theories.
Then $T$ trace defines $\Sa M$ if $\Sa M$ is trace definable in some model of $T$, $T$ trace defines $T^*$ if every model of $T^*$ is trace definable in $T$, and $T$ is \textbf{trace equivalent} to $T^*$ if $T$ trace defines $T^*$ and vice versa.
Finally, we say that two structures are trace equivalent when their theories are.
Lemma~\ref{lem:trace-transitive} is immediate.

\begin{lemma}\label{lem:trace-transitive}
Trace definability gives a transitive relation between structures and between theories.
Trace equivalence is an equivalence relation between structures and between theories.
\end{lemma}

We often use Lemma~\ref{lem:trace-transitive} without reference.
Recall that two structures or theories are {\bf mutually interpretable} if each interprets the other.




\begin{proposition}
\label{prop:trace-interpret}
If $\Sa N$ interprets $\Sa M$ then $\Sa N$ trace defines $\Sa M$.
If $T$ interprets $T^*$ then $T$ trace defines $T^*$.
Mutually interpretable structures or theories are trace equivalent.
\end{proposition}


\begin{proof}
We prove the first claim.
Suppose $\Sa N$ interprets $\Sa M$.
Let $X \subseteq N^n$ be $\Sa N$-definable, $E\subseteq X^2$ be an $\Sa N$-definable equivalence relation, and $\uppi \colon X \to M$ be a surjection so that
\begin{enumerate} 
\item for all $\alpha,\beta \in X$ we have $E(\alpha,\beta)$ if and only if $\uppi(\alpha) = \uppi(\beta)$, and
\item the pullback of any $\Sa M$-definable set by $\uppi$ is $\Sa N$-definable.
\end{enumerate}
Let $\uptau$ be a section of $\uppi$.
If $Y$ is an $\Sa M$-definable set and  $Y^*$ is the pullback of $Y$ by $\uppi$ then $Y^*$ is $\Sa N$-definable and $Y$ is the pullback of $Y^*$ by $\uptau$.
So $\uptau$ is a trace definition $\Sa M \rightsquigarrow \Sa N$.
\end{proof}

It is sufficient to realize zero-definable sets as traces.

\begin{lemma}\label{lem:b2}
Suppose $\uptau \colon M \hookrightarrow N^n$ is such that every set that is $\Sa M$-definable without parameters is a pullback of an $\Sa N$-definable set by $\uptau$.
Then $\uptau$ is a trace definition.
\end{lemma}

The proof of Lemma~\ref{lem:b2} is easy and left to the reader.

\begin{lemma}\label{lem:ka}
If $\Sa N$ trace defines $\Sa M \models T$ then $\Sa M$ is trace definable in the reduct of $\Sa N$ to some sublanguage of cardinality $\le |T|$.
\end{lemma}

\begin{proof}
Let $\uptau$ be a trace definition $\Sa M \rightsquigarrow \Sa N$.
Let $L$ be the language of $\Sa N$.
For every set $X$ that is zero-definable in $\Sa M$ fix an $\Sa N$-definable set $Y_X$ such that $X$ is the pullback of $Y_X$ by $\uptau$.
Let $L'$ be the smallest sublanguage of $L$ such that every $Y_X$ is $L'(N)$-definable.
Then $|L'|\le |T|$ and an application of Lemma~\ref{lem:b2} shows that $\uptau$ is a trace definition $\Sa M \rightsquigarrow \Sa N \!\upharpoonright\! L'$.
\end{proof}

We now show that if some model of $T^*$ is trace definable in $T$ then every model of $T^*$ is trace definable in $T$.
(Recall that all theories are assumed to be complete.)

\begin{lemma}\label{lem:b3}
The following are equivalent for any theories $T,T^*$.
\begin{enumerate}[leftmargin=*]
\item $T$ trace defines $T^*$.
\item Some model of $T^*$ is trace definable in a model of $T$.
\end{enumerate}
\end{lemma}

It follows that two structures are trace equivalent if and only if each is trace definable in an elementary extension of the other.

\begin{proof}
It is clear that (1) implies (2).
Suppose that $\Sa M\models T^*$, $\Sa N \models T$, and $\uptau$ is a trace definition $M \rightsquigarrow \Sa N$.
Let $\Sa O$ be an arbitrary model of $T^*$.
We show that $\Sa O$ is trace definable in $T$.
Consider the two-sorted structure $(\Sa M, \Sa N, \uptau)$, and let $(\Sa M^*, \Sa N^*, \uptau^*)$ be an $|O|^+$-saturated elementary extension of $(\Sa M, \Sa N, \uptau)$.
An application of Lemma~\ref{lem:b2} shows that $\uptau^*$ is a trace definition $\Sa M^* \rightsquigarrow \Sa N^*$.
By saturation there is an elementary embedding $e\colon \Sa O \hookrightarrow \Sa M^*$.
Composing $e$ and $\uptau^*$ gives a trace definition $\Sa O\rightsquigarrow\Sa N^*$.
\end{proof}

The proof of Proposition~\ref{prop:trace-theories} shows that if there is an $m$-ary trace definition of some model of $T^*$ into a model of $T$ then there is an $m$-ary trace definition of any model of $T^*$ into a model of $T$.
We therefore say that there is an $m$-ary trace definition of $T^*$ in $T$ if there is an $m$-ary trace definition of some model of $T^*$ in a model of $T$.

\medskip
Proposition~\ref{prop:trace-theories} also follows by the proof of Lemma~\ref{lem:b3}.

\begin{proposition}\label{prop:trace-theories}
Suppose that $\lambda$ is an infinite cardinal and $\Sa M$ is a structure of cardinality $\le \lambda$ in a language of cardinality $\le \lambda$.
If $\Sa N$ is a $\lambda^+$-saturated model of a theory $T$ then $T$ trace defines $\Sa M$ if and only if $\Sa N$ trace defines $\Sa M$.
\end{proposition}

If $\Sa M$ admits quantifier elimination in a relational language then to trace define $\Sa M$ it is sufficient to realize the basic relations as traces.

\begin{proposition}
\label{prop:qe}
Suppose that $L^*$ is relational,  $\Sa M$ is an $L^*$-structure with quantifier elimination, and $\Sa N$ is an arbitrary structure.
Then the following are equivalent:
\begin{enumerate}
[leftmargin=*]
\item $\Sa N$ trace defines $\Sa M$.
\item $\Sa M$ embeds into an $\Sa N$-definable $L^*$-structure.
\item Up to isomorphism $M \subseteq N^n$ for some $n \ge 1$ and $\{ a \in M^m : \Sa M \models R(a) \}$ is a trace of an $\Sa N$-definable set for every $m$-ary $R \in L^*$. 
\item There is an injection $\uptau \colon M \hookrightarrow N^n$ such that $\{ a \in M^m : \Sa M \models R(a) \}$ is a pullback of an 
$\Sa N$-definable subset of $N^{mn}$ for every $m$-ary relation $R \in L^*$.
\end{enumerate}
\end{proposition}

\begin{proof}
We leave it to the reader to show that (2), (3), and (4) are equivalent.
Propositions~\ref{prop:qe-trace} and \ref{prop:trace-interpret} together show that (2) implies (1).
Finally, it is clear from the definition that any trace definition $\Sa M\rightsquigarrow\Sa N$ satisfies the condition imposed on $\uptau$ in (4).
\end{proof}




\begin{proposition}
\label{prop:stone}
Suppose that there is an $m$-ary trace definition of an $L^*$-theory $T^*$  in $T$.
Fix $\Sa M\models T^*$, a set $B\subseteq M$ of parameters, and $n\ge 1$.
Then $S_n(\Sa M,B)$ is a continuous image of a closed subset of $S_{nm}(\Sa N,A)$ for some $\Sa N\models T$ and set $A\subseteq N$ of parameters of cardinality at most $|B|+|T^*|$.  
\end{proposition}

\begin{proof}
After possibly replacing $\Sa M$ with $\Sa M[n]$ and $nm$ with $m$ we suppose that $n=1$.
After possibly adding constant symbols for the elements of $A$ to $L^*$, and noting that the resulting structure is still trace definable in $T$, we suppose that $B=\varnothing$.
We may also suppose that $M\subseteq N^m$ and the inclusion $M \hookrightarrow N^m$ is a trace definition.
Let $P$ be an $m$-ary relation on $N$ defining $M$.
Let $\lambda= |T| + |T^*|$.
By the proof of Lemma~\ref{lem:b3} we may suppose that $(\Sa N,P)$ is $\lambda^+$-saturated.
In particular $\Sa M$ is $|T^*|^+$-saturated.
Fix a set $A\subseteq N$ of parameters such that $|A|=|T^*|$ and every zero-definable subset of every $M^k$ is a trace of an $A$-definable subset of $N^{mk}$.
It follows that if $\beta,\beta^*\in M$ then $\tp_{\Sa N}(\beta|A)=\tp_{\Sa N}(\beta^*|A)$ implies $\tp_{\Sa M}(\beta)=\tp_{\Sa M}(\beta^*)$.
Let $Y$ be $\{ \tp_{\Sa N}(\beta|A) : \beta\in M\}$ and let $f\colon Y \to S_1(\Sa M)$ be the map given by declaring $f(p)=q$ when there is $\beta \in M$ such that $p=\tp_{\Sa N}(\beta|A)$ and $q=\tp_{\Sa M}(\beta)$.
Then $Y$ is closed by $\lambda^+$-saturation of $(\Sa N,P)$ and $f$ is surjective as $\Sa M$ is $\lambda^+$-saturated.
We show that $f$ is continuous.
Let $U\subseteq S_1(\Sa M)$ be clopen.
Then $U$ is the set of types which concentrate on some zero-definable $X\subseteq M$.
Fix an $A$-definable set $Z\subseteq N^{m}$ such that $X=Z\cap M$.
Then $f^{-1}(U)$ is the set of types $p \in Y$ which concentrate on $Z$.
Hence $f^{-1}(U)$ is clopen in $Y$.
\end{proof}

\subsection{Embeddings between Tarski systems}\label{section:ts}
This section is not needed for the rest of the paper.
We show that trace embeddability is equivalent to embeddability between ``language free structures" in a precise sense.
We first recall the classical notion of a ``language free embedding".
This would usually  be phrased in terms of cylindrical set algebras, see \cite{Monk1993}.
We state it in a slightly different form to match our convention that ``definable" means ``definable, possibly with parameters".

\medskip
A \textbf{Tarski system} $\Cal S$ is a set $S$ together with a family $\Cal S_n$ of subsets of $S^n$ for each $n\ge 1$ such that the following are satisfied for all $1 \le n < m$.
\begin{enumerate}[leftmargin=*]
\item $\Cal S_n$ is a boolean algebra of subsets of $S^n$.
\item $\{a\}\in\Cal S_1$ for each $a\in S$.
\item $\{(a_1, a_2)\in S^2 : a_1 = a_2\} \in \Cal S_2$.
\item If $\uppi$ is a coordinate projection $S^m \to S^{n}$ then $\uppi(X) \in \Cal S_n$ when $X \in \Cal S_m$ and $\uppi^{-1}(X) \in \Cal S_m$ when $X \in \Cal S_n$.
\end{enumerate}
We associate a Tarski system $\Cal S (\Sa M)$ to a structure $\Sa M$ in the obvious way.
Any Tarski system arises in this way.
We say that an algebraic embedding  $\uptau \colon\Cal S\hookrightarrow\Cal U$ of Tarski systems is an injection  $\uptau\colon S\hookrightarrow U$ together with an injection $\uptau_n\colon \Cal S_n\hookrightarrow\Cal U_n$ for each $n\ge 1$ so that:
\begin{enumerate}[leftmargin=*, label=(\alph*)]
\item Each $\uptau_n$ is an embedding of boolean algebras.
\item If $X \in \Cal S_{n}$, $i \in \{1, \ldots, n\}$, and $\uppi, \uppi^*$ are the projections $S^n \to S^{n - 1}$, $U^n \to U^{n - 1}$ away from the $i$th coordinate, respectively, then $\uppi^*(\uptau_n(X)) = \uptau_{n - 1}(\uppi(X))$.
\item We have the following for all $X\in\Cal S_n$ and $\alpha_1,\ldots,\alpha_n\in S$:
$$(\alpha_1,\ldots,\alpha_n)\in X\quad\Longleftrightarrow\quad (\uptau(\alpha_1),\ldots,\uptau(\alpha_n))\in \uptau_n(X)$$
\end{enumerate}

\begin{fact}
\label{fact;ee}
The following are equivalent for any structures $\Sa M, \Sa N$ and $\uptau \colon M \to N$.
\begin{enumerate}[leftmargin=*]
\item $\uptau$ extends to an algebraic embedding $\Cal S(\Sa M)\hookrightarrow\Cal S(\Sa N)$.
\item $\uptau$ gives an elementary embedding of $\Sa M$  into a reduct of $\Sa N$.
\end{enumerate}
\end{fact}

See \cite[Thm~7.1]{Monk1993} for the classical version of this fact.

\begin{proof}[Proof sketch]
Suppose (2).
We may suppose that $\Sa M$ is an elementary substructure of a reduct $\Sa N^*$ of $\Sa N$ and $\uptau$ is the inclusion.
For each $\Sa M$-definable subset $X\subseteq M^n$ let $\uptau_n(X)$ be the subset of $N^n$ defined in $\Sa N^*$ by any formula that defines $X$.
It is easy to see that this gives an algebraic embedding $\Cal S(\Sa M) \hookrightarrow \Cal S(\Sa N)$.

\medskip
Suppose (1) and let $(\uptau_n)_{n \ge 1}$ be the extension.
We may suppose that $M \subseteq N$ and $\uptau$ is the inclusion.
Then $X = \uptau_n(X)\cap M^n$ for every $\Sa M$-definable $X\subseteq M^n$ by (c).
After possibly Morleyizing we suppose that $\Sa M$ admits quantifier elimination and the language $L$ of $\Sa M$ is relational.
Let $\Sa N^*$ be the $L$-structure on $N$ given by declaring that $a \in N^n$ satisfies an $n$-ary $R \in L$ if and only if $a$ lies in the set given by applying $\uptau_n$ to the subset of $M^n$ defined by $R$.
Then $\Sa N^*$ is a reduct of $\Sa N$ and $\Sa M$ is a substructure of $\Sa N^*$.
The Tarski-Vaught test  shows that $\Sa M$ is an elementary substructure of $\Sa N^*$.
It suffices fix a quantifier-free $L(m)$-formula $\phi(x_1, \ldots, x_n)$ and $a \in M^{n - 1}$ and show that  $\Sa M \models \exists x  \phi(a, x)$ if and only if $\Sa N^* \models \exists x \phi (a, x)$.
This follows by applying (b) with $X$ the set of $b \in M^{n}$ such that $\Sa M \models \exists x \phi(b)$.
\end{proof}

We now define a ``non-algebraic" notion of embedding.
An embedding $\uptau\colon\Cal S\hookrightarrow\Cal U$ of Tarski systems consists of an injection $\uptau\colon S\hookrightarrow U$ together with an injection $\uptau_n\colon \Cal S_n\hookrightarrow\Cal U_n$ for each $n\ge 1$ such that (c) above holds.
The assumption of injectivity on the $\uptau_n$  is superfluous as (c) ensures that every $X\in\Cal S_n$ is determined by $\uptau_n(X)$.


\begin{proposition}\label{prop:tarski system}
The following are equivalent for any structures $\Sa M, \Sa N$ and  $\uptau \colon M \to N$.
\begin{enumerate}[leftmargin=*]
\item $\uptau$ extends to a Tarski system embedding $\Cal S(\Sa M) \hookrightarrow \Cal S(\Sa N)$.
\item $\uptau$ is a trace embedding $\Sa M \hookrightarrow \Sa N$.
\end{enumerate}
\end{proposition}

\begin{proof}
It is clear from the definitions that (1) implies (2).
Suppose (2).
We extend $\uptau\colon\Sa M\hookrightarrow\Sa N$  to a Tarski system embedding $\Cal S(\Sa M)\hookrightarrow\Cal S(\Sa N)$ by letting $\uptau_n(X)$ be an arbitrary $\Sa N$-definable subset of $N^n$ whose pullback under $\uptau$ is $X$ for every $n\ge 1$ and $\Sa M$-definable $X\subseteq M^n$.
\end{proof}

Finally, there is an $\aleph_0$-sorted language $L$ and an incomplete $L$-theory $T$ such that the category of  Tarski system embeddings is equivalent to the category of embeddings between models of $T$.
We leave the details of this to the reader.

\subsection{Unary structures}\label{section:unar}
Recall that a theory is unary if every formula is equivalent to a boolean combination of unary formulas and formulas in the language of equality.
For convenience we suppose that all unary structures are infinite.
A \textbf{unary relational structure} is a structure in a unary relational language.
It is clear that any unary structure is interdefinable with a unary relational structure.

\begin{fact}
\label{fact:unary}
Any unary relational structure is homogeneous and admits quantifier elimination.
\end{fact}

\begin{proof}
It is enough to prove the first claim.
Let $\Sa X$ be a unary $L$-structure, $a_1,a^*_1,\ldots,a_n,a^*_n$ be elements of $X$ such that the $a_i$ are distinct and the $a^*_i$ are also distinct, and suppose that we have $R(a_i)$ if and only if $R(a^*_i)$ for all unary $R\in L$ and $i=1,\ldots,n$.
Let $\sigma\colon X \to X$ be given by $\sigma(a_i)=a^*_i$ and $\sigma(a^*_i)=a_i$ for  $1 \le i \le n$ and $\sigma(b)=b$ when $b\notin \{a_1,a^*_1,\ldots,a_n,a^*_n\}$.
Note that $\sigma$ is an automorphism of $\Sa X$.
\end{proof}

Fact~\ref{fact:unary}  has the following consequences:
\begin{enumerate}
[leftmargin=*]
\item Any unary relational structure is indeed a unary structure.
\item A structure is unary if and only if it is interdefinable with a unary relational structure.
\end{enumerate}

\begin{lemma}\label{lem:un-trace}
Let $L$ be unary relational, $\Sa X$ be an $L$-structure, and $\Sa M$ be a structure.
Then $\uptau\colon X \to M^n$ is a trace definition if and only if $\uptau$ is injective and $\{ a \in X : \Sa X \models R(a) \}$ is the pullback of an $\Sa M$-definable subset of $M^n$ for every unary $R \in L$.
\end{lemma}

Lemma~\ref{lem:un-trace} follows from Proposition~\ref{prop:qe-trace} and quantifier elimination for unary relational structures.
We now show in particular that a unary structure is determined up to trace equivalence by its boolean algebra of zero-definable unary sets.

\begin{lemma}\label{lem:both ways}
Let $\Sa X, \Sa X^*$ be unary structures and suppose that there is either a continuous injection $S_1(\Sa X) \hookrightarrow S_1(\Sa X^*)$ or a continuous surjection $S_1(\Sa X^*) \to S_1(\Sa X)$.
Then $\Th(\Sa X^*)$ trace defines $\Th(\Sa X)$.
Hence $\Sa X$ is trace equivalent to $\Sa X^*$ when $S_1(\Sa X)$ is homeomorphic to $S_1(\Sa X^*)$.
\end{lemma}

\begin{proof}
We first reduce to the case when every $p \in S_1(\Sa X^*)$ has at least $|X|$ realizations in $X^*$.
After possibly replacing $\Sa X^*$ with an elementary extension we may suppose that $\Sa X^*$ is $\max(|X|, |L^*|)^+$-saturated.
Let $\Sa X^{**}$ be the $L^*$-structure with domain $X^* \times X^*$ where for every unary $R \in L^*$ and $(a, b) \in X^* \times X^*$ we have $\Sa X^{**} \models R(a,b)$ if and only if $\Sa X^* \models R(a)$.
Note that every $p \in S_1(\Sa X^{**})$ has at least $|X|^+$-realizations in $\Sa X^{**}$.
It is clear that $\Sa X^{**}$ is interpretable in $\Sa X^*$.
Furthermore if $b$ is an arbitrary element of $X^*$ then the map $a \mapsto (b, a)$ is an embedding $\Sa X^* \hookrightarrow \Sa X^{**}$ and is hence a trace embedding by quantifier elimination for unary relational structures.
Hence $\Sa X^{**}$ is trace equivalent to $\Sa X^*$.
So after possibly replacing $\Sa X^*$ with $\Sa X^{**}$ we have the desired reduction.

\medskip
Suppose that $f$ is a continuous injection $S_1(\Sa X) \hookrightarrow S_1(\Sa X^*)$.
Let $\uptau$ be an injection $X \hookrightarrow X^*$ such that $\uptau(a)$ realizes $f(\tp_{\Sa X}(a))$ for all $a \in X$.
We show that $\uptau$ is a trace embedding $\Sa X \hookrightarrow \Sa X^*$ by applying Lemma~\ref{lem:un-trace}.
Fix unary $R \in L$.
Let $U$ be the set of $p(x) \in S_1(\Sa X)$ containing $R(x)$, so $U$ is clopen.
Now $f(U)$ is a clopen subset of $f(S_1(\Sa X))$, so there is clopen $V \subseteq S_1(\Sa X^*)$ such that $f^{-1}(V) = U$.
Let $\phi(x)$ be an $L^*$-formula such that $V$ is the set of $p(x) \in S_1(\Sa X^*)$ containing $\phi(x)$.
Tracing through, we have $\Sa X \models R(a)$ if and only if $\Sa X^* \models \phi(\uptau(a))$ for all $a \in X$.
Hence $\uptau$ is a trace embedding.

\medskip
Now suppose that $f$ is a continuous surjection $S_1(\Sa X^*) \to S_1(\Sa X)$.
Let $\uptau$ be an injection $X \to X^*$ such that $f(\tp_{\Sa X^*}(b)) = \tp_{\Sa X}(a)$ when $b = \uptau(a)$.
We again apply Lemma~\ref{lem:un-trace} to show that $\uptau$ is a trace embedding $\Sa X \hookrightarrow \Sa X^*$.
Fix unary $R \in L$ and let $U \subseteq S_1(\Sa X)$ be as in the previous paragraph.
Let $V = f^{-1}(U)$.
Note that $V$ is a clopen subset of $S_1(\Sa X^*)$ and let $\phi(x)$ be as in the previous paragraph.
Now we have $\Sa X \models R(a)$ if and only if $\Sa X^* \models \phi(\uptau(a))$ for all $a \in X$.
Hence $\uptau$ is a trace embedding.
\end{proof}

Proposition~\ref{prop:unary char} shows that the class of unary structures trace definable in a theory $T$ is determined solely by the topologies of type spaces over $T$.

\begin{proposition}
\label{prop:unary char}
The following are equivalent for an arbitrary theory $T$ and unary $L^*$-structure $\Sa X$.
\begin{enumerate}
\item $T$ trace defines $\Sa X$.
\item $S_1(\Sa X)$ is a continuous image of a closed subset of $S_n(\Sa M, A)$ for some $\Sa M\models T$, $A\subseteq M$, and $n\ge 1$.
\item  $S_1(\Sa X)$ is a continuous image of a closed subset $Y$ of $S_n(\Sa M, A)$ for some $\Sa M\models T$, $A\subseteq M$ such that $|A|\le |T^*|$, and $n\ge 1$.
\end{enumerate}
\end{proposition}


\begin{proof}
Proposition~\ref{prop:stone} shows that (1) implies (3).
It is clear that (3) implies (2).
Suppose (2).
Let $\Sa M^*$ be the unary relational structure with domain $M^n$ and relations defining all $A$-definable  subsets of $M^n$.
Then $\Sa M^*$ is interpretable in $\Sa M$ and $S_1(\Sa M^*) = S_n(\Sa M, A)$.
Let $\Sa Y$ be the unary relational structure with domain $Y$ and relations defining all clopen subsets of $Y$.
Then there is a continuous injection $S_1(\Sa Y) \hookrightarrow S_1(\Sa M^*)$, so $\Th(\Sa M^*)$ trace defines $\Sa Y$ by Lemma~\ref{lem:both ways}.
Finally, there is a continuous surjection $S_1(\Sa Y) \to S_1(\Sa X)$, so $\Th(\Sa Y)$ trace defines $\Sa X$ by another application of Lemma~\ref{lem:both ways}.
\end{proof}

\section{Stability, $k$-$\nip$, and indiscernibles}\label{section:sick}
We first characterize stability.
We say that a unary structure $\Sa X$ is \textbf{chainable} if $\Sa X$ is interdefinable with $(X;\Cal G)$ for some collection $\Cal G$ of subsets of $X$ which forms a chain under inclusion.
Note that $\Sa X$ is chainable if and only if the boolean algebra of unary definable sets is generated by a chain.
This is a well-known condition in the theory of boolean algebras.
The following are equivalent for any boolean algebra $\mathfrak{B}$.
\begin{enumerate}[leftmargin=*]
\item $\mathfrak{B}$ is generated by a chain.
\item $\mathfrak{B}$ is isomorphic to the boolean algebra generated by all intervals in some linear order.
\item The topology on the Stone space $S(\mathfrak{B})$ of $\mathfrak{B}$ is induced by a linear order $\triangleleft$, and in this case $\mathfrak{B}$ is generated by the collection of intervals of the form $(-\infty,\upxi]$ for $\upxi\in S(\mathfrak{B})$.
\item There is a linear order $(I;\triangleleft)$ such that $\mathfrak{B}$ is isomorphic to the boolean algebra of subsets of $I$ generated by all sets of the form $(-\infty,\upxi]$ for $\upxi\in I$.
\end{enumerate}
The equivalence of (1) and (2) is \cite[Thm~15.3]{booleanhandbook}.
The equivalence of (2) and (3) follows by \cite[Thm~15.7]{booleanhandbook} and the proof of that theorem.
It is easy to see that (3) implies (4) and that (4) implies (1).
Given a linear order $(I;\triangleleft)$ let $\Sa I_\triangleleft$ be the unary relational structure with domain $I$ and unary relations defining $(-\infty,\alpha]$ for all $\alpha\in I$.
Lemma~\ref{lem:chain0} follows by Lemma~\ref{lem:both ways} and  the comments above.

\begin{lemma}
\label{lem:chain0}
Any chainable unary structure is trace equivalent to $\Sa I_\triangleleft$ for a linear order $(I;\triangleleft)$.
\end{lemma}

We say that a unary structure $\Sa X$ is {\bf treeable} if $\Sa X$ is interdefinable with $(X; \Cal G)$ for some collection $\Cal G$ of subsets which form a tree under inclusion, i.e. if two elements of $\Cal G$ intersect then one is contained in the other.
We now characterize stability.

\begin{proposition}
\label{prop:stable-0}
Stability is preserved under trace definability.
Furthermore the following are equivalent for an arbitrary theory $T$.
\begin{enumerate}[leftmargin=*]
\item $T$ is unstable
\item $T$ trace defines $(\Q;<)$.
\item $T$ trace defines an infinite linear order.
\item $T$ trace defines every chainable unary structure.
\item $T$ trace defines $\Sa I_\triangleleft$ for some linear order $(I; \triangleleft)$ such that $|T|<|I|$.
\item $T$ trace defines every treeable unary structure.
\end{enumerate}
\end{proposition}

In particular it follows that a countable theory is unstable if and only if it trace defines the unary relational structure with domain $\omega_1$ and relations defining $[0, \eta]$ for all $\eta < \omega_1$.

\begin{proof}
We first show that (2) implies (6).
It suffices to show that $\dlo$ trace defines any treeable unary structure.
Let $\Sa X = (X; \Cal G)$ where $\Cal G$ is a collection of subsets of $X$ such that $Y \cap Y^* \in \{Y, Y^*, \varnothing\}$ for any $Y, Y^* \in \Cal G$.
By a theorem of Heindorf there is a linear order $(I ; \triangleleft)$ such that $\Cal B[\Sa X]$ embeds into the boolean algebra generated by intervals in $(I; \triangleleft)$~\cite{tr}.
Let $(J ; \triangleleft)$ be any model of $\dlo$ extending $(I; \triangleleft)$.
Then $\Cal B[\Sa X]$ embeds into the boolean algebra generated by intervals in $J$, hence $(J ; \triangleleft)$ trace defines $\Sa X$ by Proposition~\ref{prop:unary char}.

\medskip
It is clear that (6) implies (4) and that (4) implies (5).
We show that (5) implies (1).
Suppose that $\Sa I_\triangleleft$ is as in (5) and that $\Sa M\models T$ trace defines $\Sa I_\triangleleft$.
We may suppose that $I\subseteq M^n$ and that the inclusion $I \hookrightarrow M^n$ is a trace definition.
After possibly replacing $\Sa M$ with $\Sa M[n]$ we suppose that $n = 1$.
For each $\beta\in I$ fix $\Sa M$-definable $Y_\beta\subseteq M$ such that $(-\infty,\beta]=Y_\beta\cap I$.
As $|T|<|I|$ there is a parameter-free formula $\varphi(x,y)$ with $|x|=1$, a sequence $(\beta_i)_{i<\omega}$ of distinct elements of $I$, and a sequence $(\gamma_i)_{i<\omega}$ of elements of $M^{|y|}$ such that $Y_{\beta_i}=\varphi(M,\gamma_i)$ for all $i<\omega$.
By Ramsey's theorem we may suppose that $(\beta_i)_{i<\omega}$ is either strictly increasing or strictly decreasing.
Without loss of generality we only treat the strictly increasing case.
Then for all $i,j\in I$ we have $\beta_i\in Y_{\beta_j}$ if and only if $i\triangleleft j$.
Hence we have $\Sa M \models \varphi(\beta_i,\gamma_j)$ if and only if $i\triangleleft j$.
Therefore $\varphi$ is unstable.

\medskip
We show that (1) implies (2).
By Proposition~\ref{prop:trace-theories} it is enough to suppose that $\Sa M$ is unstable and $\aleph_1$-saturated and show that $\Sa M$ trace defines $(\Q;<)$.
There is a tuple $(\alpha_q)_{q \in \Q}$ of elements of some $M^m$ and a formula $\phi(x,y)$ such that for all $p,q \in \Q$ we have $\Sa M \models \phi(\alpha_p,\alpha_q)$ if and only if $p < q$.
Let $\uptau \colon \Q \hookrightarrow M^m$ be given by declaring $\uptau(q) = \alpha_q$ for all $q \in \Q$.
As $(\Q;<)$ admits quantifier elimination an application of Proposition~\ref{prop:qe} shows that  $\uptau$ is a trace definition $(\Q;<) \rightsquigarrow \Sa M$.

\medskip
Clearly (2) implies (3).
Finally, (3) implies (2) as infinite linear orders are unstable.
\end{proof}

We next consider $k$-$\nip$.
We first state an easy lemma.
Recall that a structure or theory is $k$-ary if every formula is equivalent to a boolean combination of formulas in the language of equality and formulas in at most $k$ variables.

\begin{lemma}\label{lem:k-ary}
Suppose that $\Sa M$ is $k$-ary for $k \ge 1$.
Then $\Sa M$ is interdefinable with an $L^*$-structure $\Sa M^*$ such that every relation in $L^*$ other than equality has arity $k$ and $\Sa M^*$ admits quantifier elimination.
Now suppose that $\Sa M$ is an $L$-structure with quantifier elimination for a finite relational language $L$ such that every relation in $L$ other than equality has arity $\le k$.
Then $\Sa M$ is interdefinable with an $L^*$-structure structure $\Sa M^*$ such that $L^*$ is finite, every relation in $L^*$ other than equality has arity exactly $k$, and $\Sa M^*$ admits quantifier elimination.
\end{lemma}

We leave the proof of Lemma~\ref{lem:k-ary} to the reader.
The idea is to Morleyize and then replace each relation $R(x_1,\ldots,x_m)$ with $m < k$ with a relation $R^*(x_1,\ldots,x_k)$ given by declaring $R^*(a_1,\ldots,a_k)$ if and only if $R(a_1,\ldots,a_m)$ and $a_m = a_{m+1} = \ldots = a_k$.

\begin{lemma}
\label{lem:random}
Fix $k\ge 2$.
Any structure admitting quantifier elimination in a finite $k$-ary relational language is trace definable in the theory of the generic $k$-ary relation.
\end{lemma}

\begin{proof}
Suppose $\Sa M=(M;R_1,\ldots,R_n)$ admits quantifier elimination and each $R_i$ has arity $\le k$.
We show that $\Sa M$ is trace definable in the theory of the generic $k$-ary relation.
We may suppose that $\Sa M$ is countable by Lemma~\ref{lem:b3}.
By Lemma~\ref{lem:k-ary} we may suppose that every $R_i$ has arity exactly $k$.
Let $(N;R)$ be the generic $k$-ary relation.
Then any countable $k$-ary relation embeds into $(N;R)$ so we fix an embedding $\uptau_i\colon (M;R_i) \hookrightarrow (N;R)$ for each $i = 1,\ldots,n$.
By Proposition~\ref{prop:qe} the map $\uptau \colon M \hookrightarrow N^n$ given by $\uptau = (\uptau_1,\ldots,\uptau_n)$ is a trace definition $\Sa M \rightsquigarrow (N;R)$.
\end{proof}

We now characterize $k$-$\nip$.
See \cite{cpt} for background on $k$-$\nip$ and recall that $1$-$\nip$ is $\nip$.

\begin{proposition}\label{prop:k nip}
Fix $k \ge 1$.
The class of $k$-$\nip$ structures is closed under trace definability.
Furthermore the following are equivalent for any theory $T$.
\begin{enumerate}[leftmargin=*]
\item $T$ is $k$-$\ip$.
\item $T$ trace defines the generic $(k+1)$-hypergraph.
\item $T$ trace defines the generic ordered $(k+1)$-hypergraph.
\item $T$ trace defines the generic $(k+1)$-ary relation.
\item $T$ trace defines any structure that admits quantifier elimination in a finite $(k+1)$-ary relational language.
\item $T$ trace defines any $k$-ary structure.
\item $T$ trace defines the model companion of the theory of a set equipped with $|T|^+$ relations, each of arity $k$.
\end{enumerate}
\end{proposition}


\begin{proof}
Let $d = k + 1$ and $\kappa = |T|^+$.
Preservation of $k$-$\nip$ under trace definability follows easily from the definition of $k$-$\nip$.
Hence each of (2), (3), (4) implies (1) as each of these structures is $k$-$\ip$.
Furthermore (5) implies each of (2), (3), (4) as each structure admits quantifier elimination in a finite $k$-ary relational language.
Lemma~\ref{lem:random} shows that (4) implies (5).
Let $T^*$ be the model companion described in (7).
Then $T^*$ exists, is complete, and admits quantifier elimination by~\cite[Thm.~3.7, Cor.~B.3]{jera}.
In particular $T^*$ is a $k$-ary theory, so (6) implies (7).

\medskip
We show that (7) implies (4).
Suppose that (7) holds.
Fix $\Sa M \models T$ and $\Sa O = (O; (R_i)_{i < \kappa}) \models T^*$ and suppose that $\Sa M$ trace defines $\Sa O$.
After possibly replacing $\Sa M$ with $\Sa M[n]$ for some $n$ we may suppose that $\Sa O$ trace embeds into $\Sa M$.
Hence we may suppose that $O$ is a subset of $M$ and the inclusion $\Sa O \to \Sa M$ is a trace embedding.
For every $i < \kappa$ fix a formula $\varphi_i(x_1,\ldots,x_{n_i},y_1,\ldots,y_k)$, and $\beta_i \in M^{n_i}$ such that we have $\Sa M \models \varphi_i(\beta_i, a)$ if and only if $\Sa O \models R_i(a)$ for any $a \in O^k$.
Now as $\kappa > |T| \ge \aleph_0$ there is a countably infinite $I \subseteq \kappa$ such that $i \mapsto n_i$ and $i \mapsto \varphi_i$ are constant on $I$.
Set $n = n_i$ and $\varphi = \varphi_i$ for any $i \in I$.
Let $\Sa P = (P;(R_i)_{i \in I})$ be a countable elementary  submodel of $(O;(R_i)_{i \in I})$.
Then $\Sa P$ is a model of the model companion of the theory of a set equipped with $\aleph_0$ relations, each of arity $k$.
Fix a bijection $\sigma \colon P \to I$ and let $R$ be the $d$-ary relation on $P$ given by $R(a_1,\ldots,a_{k}, b)$ if and only if $R_{\sigma(b)}(a_1,\ldots,a_{k})$.
It is easy to see that any countable $k$-ary relation embeds into $(P; R)$.
Hence we there is $P^* \subseteq P$ such that $(P^*; R)$ is an isomorphic copy the generic $k$-ary relation.
Let $\iota \colon P \to M^n$ be given by $\iota(b) = \beta_{\sigma(b)}$.
So we have $R(a_1,\ldots,a_{k}, b)$ if and only if $\Sa M \models \varphi(\iota(b), a_1,\ldots,a_k)$.
It follows by an application of Proposition~\ref{prop:qe} that the map $P^* \hookrightarrow M \times M^n$ given by $b \mapsto (b, \iota(b))$ is a trace definition $(P^*; R) \rightsquigarrow \Sa M$.

\medskip
We finish by showing that (1) implies (6).
Suppose that $T$ is $k$-$\ip$.
Let $\Sa O$ be a $k$-ary structure.
By Lemma~\ref{lem:k-ary} we may suppose that $\Sa O$ admits quantifier elimination in a relational language $L$ such that every relation in $L$ other than equality has arity exactly $k$.
Let $(R_i)_{i < \kappa}$ be an enumeration of the relations in $L$ other than equality.
By Lemma~\ref{lem:b3} it suffices to show that $T$ trace defines any structure elementarily equivalent to $\Sa O$, so we may suppose that $|O| = \kappa$.
Suppose that $\Sa M \models T$ is $\kappa^+$-saturated and let $\varphi(x,y_1,\ldots,y_k)$ be a $k$-$\ip$ formula in $\Sa M$.
By \cite[Remark~4]{simon2021notestabilitynipvariable} we may suppose that $|x| = |y_1| = \cdots = |y_k| = 1$.
For each $i = 1,\ldots,k$ fix a subset $A_i$ of $M$ with cardinality $\kappa$ such that for every $I \subseteq A_1\times\cdots\times A_k$ there is $\beta \in M$ satisfying 
$$\Sa M \models \varphi(\beta,a_1,\ldots,a_k) \quad\Longleftrightarrow\quad (a_1,\ldots,a_k) \in I \quad \text{for all } a_1 \in A_1,\ldots, a_k \in A_k.$$
Fix a bijection $\uptau_i \colon O \to A_i$ for each $i = 1,\ldots,k$.
For each $i < \kappa$ fix $\beta_i \in M$ such that we have $\Sa M \models \varphi(\beta_i, \uptau_1(a_1),\ldots,\uptau_k(a_k))$ if and only if $\Sa O \models R_i(a_1,\ldots,a_k)$ for all $a_1,\ldots,a_k \in O$.
Let $\uptau_0\colon O \hookrightarrow M$ be an arbitrary injection.
Finally, an application of Proposition~\ref{prop:qe} shows that the map $O \hookrightarrow M^{k + 1}$ given by $a\mapsto (\uptau_0(a),\uptau_1(a),\ldots,\uptau_k(a))$ is a trace definition $\Sa O \rightsquigarrow \Sa M$.
(Note that $\uptau_0$ handles formulas in the language of equality.)
\end{proof}

\subsection{Indiscernible collapse}\label{section:collapse}
Are there other natural categories of ``language free morphisms" between structures in arbitrary languages?
Interpretations can be seen in this light; an interpretation of $\Sa M$ in $\Sa N$ is essentially a surjection $N^n \to M$ such that the pullback of any $\Sa M$-definable set is definable in $\Sa N$.
Along these lines one can consider the category of maps $\Sa M \to \Sa N$ such that the pullback of any $\Sa N$-definable set is $\Sa M$-definable.
It is perhaps more natural to consider the category defined the same way, but instead requiring that the pullback of any zero-definable set is zero-definable.
In particular a map $\uptau \colon \Q \to M$ gives a morphism $(\Q; <) \to \Sa M$ in the latter category if and only if $(\uptau(q))_{q \in \Q}$ is an indiscernible sequence of elements of $M$ (over the empty set).
In this section we consider the relationship between indiscernibles and trace definability.

\medskip
Let $\Sa M$ be a highly saturated structure, $A$ be a small set of parameters from $\Sa M$, $\Sa I$ be a small structure, and $\upgamma$ be a map $I \to M^m$ for some $m$.
We call $\upgamma$ a \textbf{picture} of $\Sa I$ in $\Sa M$.
Then $\upgamma$  is {\bf indiscernible} over $A$ if
$$ \tp_{\Sa I}(a) = \tp_{\Sa I}(b) \quad \Longrightarrow \quad \tp_{\Sa M}(\upgamma(a)|A) = \tp_{\Sa M}(\upgamma(b)|A)$$
for any $k \ge 1$ and $a,b \in I^k$.
Furthermore $\upgamma$  is {\bf uncollapsed} over $A$ if
$$ \tp_{\Sa I}(a) = \tp_{\Sa I}(b) \quad \Longleftarrow \quad \tp_{\monster}(\upgamma(a)|A) = \tp_{\monster}(\upgamma(b)|A)$$
for any $k \ge 1$ and $a,b \in I^k$.

\begin{lemma}\label{lem:uncom}
Let $\Sa M$, $\Sa I$, $\upgamma$ be as above and suppose that $\Sa I$ is $\aleph_0$-categorical.
\begin{enumerate}[leftmargin=*]
\item If $\upgamma$ is indiscernible over $A$ then the pullback of any $A$-definable subset of $M^{nm}$ by $\upgamma$ is zero-definable in $\Sa I$.
\item If $\upgamma$ is indiscernible and uncollapsed over $A$ then any zero-definable set in $\Sa I$ is the pullback of an $A$-definable set by $\upgamma$.
\end{enumerate}
\end{lemma}

\begin{proof}
Suppose that $\upgamma$ is indiscernible over $A$, let $Y \subseteq M^{nm}$ be $A$-definable, and $X \subseteq I^n$ be the pullback of $Y$ by $\upgamma$.
It is immediate from the definition that there is a subset $P$ of $S_n(\Sa I)$ such that $X$ is the set of $a \in I^n$ which realize some $p \in P$.
By Ryll-Nardzewski the set of realizations of each $p \in S_n(\Sa I)$ is zero-definable in $\Sa I$.

\medskip
Now suppose that $\upgamma$ is also uncollapsed over $A$.
We show that every zero-definable subset of $I^n$ is a pullback of an $A$-definable set via $\upgamma$.
Let $p_1,\ldots,p_d$ be an enumeration of $S_n(\Sa I)$.
By Ryll-Nardzewski it is enough to produce $A$-definable $Y_1,\ldots,Y_d$ such that for any $a \in I^n$ we have $\tp_{\Sa I}(\upgamma(a)) = p_i$ if and only if $\upgamma(a) \in Y_i$.
Let $a_i \in I^n$ be a realization of $p_i$ for each $i = 1, \ldots, d$.
As $\upgamma$ is uncollapsed $\tp_{\Sa M}(\upgamma(a_i)|A) \ne \tp_{\Sa M}(\upgamma(a_j)|A)$ when $i \ne j$.
Hence for $1 \le i < j \le d$ there is an $A$-definable $Y_{ij} \subseteq M^{nm}$ such that $\upgamma(a_j) \notin Y_{ij} \ni \upgamma(a_i)$.
Let $Y_i = \bigcap_{j = 1}^{d} Y_{ij}$ for each $i$.
Then $\upgamma(a_i) \in Y_i$ and $\upgamma(a_j) \notin Y_i$ for all $i \ne j$.
After possibly replacing each $Y_i$ with $Y_i \setminus \bigcup_{j \ne i} Y_j$ we suppose that the $Y_i$ are pairwise disjoint.
Fix $b \in I^n$.
We show that $\tp_{\Sa I}(b) = p_i$ if and only if $\upgamma(b) \in Y_i$.
Suppose $\tp_{\Sa I}(b) = p_i$.
As $\upgamma$ is indiscernible $\tp_{\monster}(\upgamma(b)|A) = \tp_{\monster}(\upgamma(a_i)|A)$, so $\upgamma(b) \in Y_i$.
Conversely, suppose $\upgamma(b) \in Y_i$.
Then $\upgamma(b) \notin Y_j$ when $j \ne i$, so $\tp_{\Sa I}(b) \ne p_j$ for all $j \ne i$.
Hence $\tp_{\Sa I}(b) = p_i$.
\end{proof}

\medskip
Let $L$ be a finite relational language and $\Sa I$ be an $L$-structure.
Then $\Sa I$ is {\bf Ramsey} if whenever $\Sa A, \Sa B \in \age(\Sa I)$, and the embeddings $\Sa A \hookrightarrow \Sa I$ are colored with finitely many colors, then there is an embedding $\iota \colon \Sa B \hookrightarrow \Sa I$ such that the embeddings of $\Sa A$ into the image of $\iota$ are monochromatic.
Note that Ramsey $L$-structures are closed under elementary expansions.
If $\Sa I$ is homogeneous then $\Sa I$ is Ramsey in our sense if and only if $\age(\Sa I)$ is a Ramsey class in the usual sense~\cite[Prop.~11.1.13]{Bodirsky_book}.

\begin{fact}\label{fact:ramsey}
Let $\Sa I$ be a Ramsey finitely homogeneous structure, $\Sa M$ be a structure, $\Sa J$ be an elementary extension of $\Sa I$, $\upgamma \colon J \to M^n$ be such that $(\Sa M, \Sa J, \upgamma)$ is highly saturated, and $A \subseteq M$ be small.
Then there is an embedding $e \colon \Sa I \hookrightarrow \Sa J$ such that $\upgamma \circ e$ is  indiscernible  over $A$.
\end{fact}

This is essentially due to Guingona, Hill, and Scow~\cite[Thm~2.13]{gcs}.
We give a proof to be complete.
We say that $p(x_1, \ldots, x_m) \in S_m(\Sa I)$ is separated if $p \models (x_i \ne x_j)$ when $i \ne j$.

\begin{proof}
We only treat the case $n = 1$.
Given $p \in S_m(\Sa I) = S_m(\Sa J)$ and $A$-definable $Y \subseteq M^m$ we say that an elementary substructure $\Sa O$ of $\Sa J$ is $(p, Y)$-indiscernible if whenever  $a, a^* \in J^{m}$ are realizations of $p$ then we have $\upgamma(a) \in Y$ if and only if $\upgamma(a^*) \in Y$.
Now $\Sa I$ embeds into any elementary substructure of $\Sa J$, so it suffices to produce an elementary substructure which is $(p, Y)$-indiscernible for every such $p, Y$.
Note that it suffices to produce an elementary substructure which is $(p, Y)$-indiscernible for every separated $p \in S_m(\Sa I)$ and $A$-definable $Y \subseteq M^m$.
By saturation it is enough to fix separated types $p_1 \in S_{m_1}(\Sa I), \ldots, p_d \in S_{m_d}(\Sa I)$ and $A$-definable sets $Y_1 \subseteq M^{m_1}, \ldots, Y_d \subseteq M^{m_d}$ and produce an elementary submodel $\Sa O$ of $\Sa J$ which is $(p_i, Y_i)$-indiscernible for $i = 1, \ldots, d$.
We apply induction on $d$.
Suppose that $\Sa O \prec \Sa J$ is $(p_i, Y_i)$-indiscernible for $i \le d - 1$.
By saturation of $(\Sa M, \Sa J, \upgamma)$ we may suppose that $\Sa O$ is $\aleph_1$-saturated.
Let $p = p_d$, $Y = Y_d$, and $m = m_d$.
As $p$ is separated there is $\Sa A \in \age(\Sa I)$ with domain $\{a_1, \ldots, a_m\}$ such that any $(b_1, \ldots, b_m) \in J^m$ realizes $p$ if and only if $a_i \mapsto b_i$ gives an embedding $\Sa A \hookrightarrow \Sa J$.
Color the embeddings $\iota \colon \Sa A \hookrightarrow \Sa O$ with two colors according to whether $(\upgamma(\iota(a_1)), \ldots, \upgamma(\iota(a_m))) \in Y$.
Applying the Ramsey property, we see that for any $\Sa B \in \age(\Sa I)$, there is an embedding $\iota \colon \Sa B \hookrightarrow \Sa O$ such that every embedding of $\Sa A$ into the image of $\iota$ is monochromatic.
By saturation there is an embedding $\iota \colon \Sa I \hookrightarrow \Sa J$ so that the embeddings of $\Sa A$ into the image of $\iota$ are monochromatic.
Finally, note that the image of $\iota$ is $(p_i, Y_i)$-indiscernible for each $i = 1, \ldots, d$.
\end{proof}

\begin{proposition}
\label{prop:picture}
Let $\Sa I$ be a Ramsey finitely homogeneous structure, $\Sa M$ be a highly saturated model of a theory $T$, and $m \ge 1$.
Then the following are equivalent.
\begin{enumerate}[leftmargin=*]
\item There is an $m$-ary trace definition $\Sa I \rightsquigarrow \Sa M$.
\item There is an uncollapsed indiscernible picture $I \to M^m$ of $\Sa I$ in $\Sa M$ over a small set of parameters.
\end{enumerate}
Therefore the following are equivalent.
\begin{enumerate}[leftmargin=*]
\setcounter{enumi}{2}
\item $T$ trace defines $\Sa I$.
\item There is a uncollapsed indiscernible picture of $\Sa I$ in $\Sa M$ over a small set of parameters.
\end{enumerate}
\end{proposition}

Guingona and Parnes proved essentially the same result in \cite{Guingona-Parnes}.

\begin{proof}
It is enough to show that (1) and (2) are equivalent.
Suppose (1).
By the proof of Lemma~\ref{lem:b3} there is an elementary extension $\Sa J$ of $\Sa I$ and an $m$-ary trace definition $\uptau \colon \Sa J \rightsquigarrow \Sa M$ such that the two-sorted structure $(\Sa M, \Sa J, \uptau)$ is also sufficiently saturated.
Let $A$ be a countable subset of $M$ such that every set which is zero-definable in $\Sa J$ is the pullback of an $A$-definable set via $\uptau$.
By Fact~\ref{fact:ramsey} there is an embedding $e \colon \Sa I \hookrightarrow \Sa J$ such that $\uptau^* = \uptau \circ e$ is indiscernible over $A$.
Note that every zero-definable set in $\Sa I$ is a pullback of an $A$-definable set via $\uptau^*$.
Hence $\uptau^*$ is uncollapsed over $A$.

\medskip
Now suppose that $\upgamma \colon I \to M^m$ is an uncollapsed indiscernible picture over a small set $A$ of parameters.
Lemmas~\ref{lem:b2} and \ref{lem:uncom}(2) together show that $\upgamma$ is a trace definition $\Sa I \rightsquigarrow \Sa M$.
\end{proof}

Proposition~\ref{prop:k nip},  Proposition~\ref{prop:picture}, and the Ramsey property for ordered $k$-hypergraphs together yield the following theorem of Chernikov, Palacin, and Takeuchi~\cite[Thm~5.4]{cpt}.

\begin{fact}\label{fact:cpt}
Fix $k \ge 1$.
A theory $T$ is $k$-$\nip$ if and only if any monster model of $T$ admits an uncollapsed indiscernible picture of the generic ordered $(k+1)$-hypergraph.
\end{fact}

We now consider partial collapse.
Let $\Sa I^*$ be a reduct of a structure $\Sa I$.
We say that $T$ satisfies $(\Sa I, \Sa I^*)$-collapse if whenever $\Sa M \models T$ is highly saturated and $A$ is a small set of parameters, then any $A$-indiscernible picture of $\Sa I$ in $\Sa M$ is an $A$-indiscernible picture of $\Sa I^*$.

\begin{proposition}\label{prop:partial collapse}
Suppose that $\Sa I$ is a finitely homogeneous Ramsey structure and $\Sa J^*$ is a reduct of $\Sa I$.
Then the class of theories satisfying $(\Sa I, \Sa I^*)$-collapse is closed under trace definability.
\end{proposition}

\begin{proof}
We suppose that $T^*$ is a theory without $(\Sa I, \Sa I^*)$-collapse and that $T^*$ is a theory which trace defines $T$.
We show that $T^*$ does not have $(\Sa I, \Sa I^*)$-collapse.
Let $\Sa J$ be a highly saturated elementary extension of $\Sa I$ and $\Sa J^*$ be the reduct of $\Sa J$ corresponding to $\Sa I^*$.
It follows that $T^*$ does not have $(\Sa J, \Sa J^*)$-collapse.
Suppose that $\Sa M \models T$ is highly saturated, $A$ is a small set of parameters from $M$, $\upgamma \colon J \to M^n$ is an $A$-indiscernible picture of $\Sa J$, and $\upgamma$ is not an $A$-indiscernible picture of $\Sa J^*$.
To simplify notation we only treat the case $n = 1$.
Hence there are $a, b \in J^k$ such that $\tp_{\Sa J^*}(a) = \tp_{\Sa J^*}(b)$ and $\tp_{\Sa M}(\upgamma(a) | A) \ne \tp_{\Sa M}(\upgamma(b) | A)$.
Fix an $A$-definable $Y \subseteq M^k$ such that $\upgamma(a) \notin Y \ni \upgamma(b)$.
Let $X \subseteq J^k$ be the pullback of $Y$ by $\upgamma$.
By Lemma~\ref{lem:uncom}(1) $X$ is zero-definable in $\Sa J$.
We have $a \notin X \ni b$, so $X$ is not zero-definable in $\Sa J^*$.
Now let $\uptau$ be an $m$-ary trace definition $\Sa M \rightsquigarrow \Sa N$ for $\Sa N \models T$.
We may suppose that $\Sa N$ is highly saturated.
Let $\upgamma' \colon J \to N^m$ be the composition of $\upgamma$ with $\uptau$.
Let $Z$ be a definable subset of $N^{mk}$ such that $Y$ is the pullback of $Z$ by $\uptau$.
Let $B$ be a small set of parameters from $N$ over which $Z$ is definable.
Applying Fact~\ref{fact:ramsey} let $e$ be an embedding $\Sa I \hookrightarrow \Sa J$ such that $\upgamma'' = \upgamma' \circ e$ is indiscernible over $B$.
It suffices to show that $\upgamma''$ is not an indiscernible picture of $\Sa I^*$ over $B$.
Now $X$ is the pullback of $Z$ by $\upgamma'$, hence $X \cap I^k$ is the pullback of $Z$ by $\upgamma''$.
Now $X \cap I^k$ is zero-definable in $\Sa I$, but is not zero-definable in $\Sa I^*$.
Hence there are $c, d \in I^k$ such that $\tp_{\Sa I^*}(c) = \tp_{\Sa I^*}(d)$ and $c \notin X \cap I^k \ni d$.
We have $\upgamma(c) \notin Z \ni \upgamma(d)$, so $\tp_{\Sa N}(\upgamma(c) | A) \ne \tp_{\Sa N}(\upgamma(d) | A)$.
\end{proof}

We expect that there are not so many finitely homogeneous structures modulo trace equivalence.
Recall that there are continuum many homogeneous directed graphs~\cite{Henson_1972}.

\begin{conj}\label{conj:fh}
\hspace{.00000000000000000000000000000000000000000000000000000000000000001cm}
\begin{enumerate}[leftmargin=*]
\item There are only countably many finitely homogeneous structures modulo trace equivalence.
More strongly: for any $k$ there are only finitely many structures modulo trace equivalence which admit quantifier elimination in a finite $k$-ary relational language.
\item Any binary finitely homogeneous structure is trace equivalent to either the trivial theory, $\dlo$, or the \Erdos-Rado graph.
Equivalently: if $\Sa M$ is a binary finitely homogeneous structure then the class of theories that do not trace define $\Sa M$ is either the class of theories of finite structures, the class of stable theories, or the class of $\nip$ theories.
\end{enumerate}
\end{conj}

It follows from Lemma~\ref{lem:random} and Proposition~\ref{prop:k nip} that any $\ip$ binary finitely homogeneous structure is trace equivalent to the \Erdos-Rado graph.
In previous work we showed that any stable binary finitely homogeneous structure is trace equivalent to the trivial theory.
Hence (2) holds if and only if every $\nip$ binary finitely homogeneous structure is trace definable in $\dlo$.
We also conjectured that any $\aleph_0$-stable $\aleph_0$-categorical structure is trace equivalent to the trivial structure.
This conjecture has recently been proven by Bertalan Bodor.
This implies in particular that any stable finitely homogeneous structure is trace equivalent to the trivial structure, giving evidence for Conjecture~\ref{conj:fh}(1).




\section{More stability-theoretic properties}
We next give preservation and characterization results for superstability and total transcendence.
Given a structure $\Sa M$, cardinal $\lambda$, and theory $T,$ we let $S_n(\Sa M,\lambda)$ be the supremum of $\{ |S_n(\Sa M,A)| : A \subseteq M, |A| \le \lambda \}$ and $S(T,\lambda)$ be the supremum of $\{ S_n(\Sa M,\lambda) : n \in \N, \Sa M \models T \}$.

\medskip
Proposition~\ref{prop:lambda-0} follows by Proposition~\ref{prop:stone}.

\begin{proposition}
\label{prop:lambda-0}
Suppose that $\lambda$ is an infinite cardinal, $T$ is an arbitrary theory, and $T^*$ is a theory of cardinality $\le \lambda$ trace definable in $T$.
Then $S(T^*,\lambda) \le S(T,\lambda)$.
Hence if $T$ is $\lambda$-stable then $T^*$ is $\lambda$-stable.
\end{proposition}

We now show that superstability and total transcendence are preserved.

\begin{corollary}
\label{cor:superstable}
Suppose that $T^*$ is trace definable in $T$.
\begin{enumerate}[leftmargin=*]
\item If $T$ is superstable then $T^*$ is superstable.
\item If $T^*$ is countable and $T$ is $\aleph_0$-stable then $T^*$ is $\aleph_0$-stable.
\item If $T$ is totally transcendental then $T^*$ is totally transcendental.
\end{enumerate}
\end{corollary}

\begin{proof}
(1) and (2) are immediate from Proposition~\ref{prop:lambda-0}.
We prove (3).
Let $L,L^*$ be the language of $T,T^*$, respectively.
We apply the fact that a structure is totally transcendental if and only if every reduct to a countable sublanguage is $\aleph_0$-stable.
Suppose that $\Sa N \models T$ trace defines $\Sa M \models T^*$.
Fix countable $L^{*}_0 \subseteq L^*$.
By Lemma~\ref{lem:ka} there is a countable sublanguage $L_0$  of $L$ such that $\Sa N\!\upharpoonright\! L_0$ trace defines $\Sa M \!\upharpoonright\! L^*_0$.
By (2) $\Sa M \!\upharpoonright\! L^*_0$ is $\aleph_0$-stable.
\end{proof}

We let $\Sa C$ be the unary relational structure with domain the Cantor set and relations defining all clopen subsets.
James Hanson essentially showed that a theory is totally transcendental if and only if it trace defines $\Sa C$~\cite{hanson-mo}.
We extend this.



\begin{proposition}
\label{prop:unary}
The following are equivalent for any theory $T$.
\begin{enumerate}[leftmargin=*]
\item $T$ trace defines $\Sa C$.
\item $T$ is not totally transcendental.
\item $T$ trace defines any unary structure in a countable language.
\end{enumerate}
\end{proposition}


\begin{proof}
It is clear that (3) implies (1) and Corollary~\ref{cor:superstable}(3) shows that (1) implies (2) as $\Sa C$ is not totally transcendental.
Suppose that $T$ is not totally transcendental.
We show that $T$ trace defines $\Sa C$.
There is $\Sa M \models T$ and a family $\Cal S$ of $\Sa M$-definable subsets of $M$ which form an infinite complete binary tree under inclusion.
Let $\Sa X$ be $(M; \Cal S)$.
Then $\Sa X$ is a unary relational structure and $\Cal B[\Sa X]$ is a countable atomless boolean algebra.
Hence $\Sa X$ is trace equivalent to $\Sa C$ by Lemma~\ref{lem:both ways}.
It remains to show that (1) implies (3).
Let $\Sa X$ be a unary structure in a countable language.
Then $\Cal B[\Sa X]$ is a countable boolean algebra, so there is an embedding $\Cal B[\Sa X] \hookrightarrow \Cal B[\Sa C]$.
Hence $\Th(\Sa C)$ trace defines $\Th(\Sa X)$ by Lemma~\ref{lem:both ways}.
\end{proof}

We now characterize superstability.
We first recall some background.
Let $\ru$ be the U-rank.
Given a structure $\Sa M$, a {\bf global type} over $\Sa M$ is a type with parameters from $M$.

\begin{fact}
\label{fact:U rank}
Suppose that $T$ is superstable and $\Sa M\models T$ is $\kappa$-saturated for an infinite cardinal $\kappa$.
For each ordinal $\lambda$ we let $S_\lambda$ be the set of global one-types over $\Sa M$ with U-rank $\ge \lambda$.
\begin{enumerate}
\item If $\ru(\Sa M)>\lambda$ then $|S_\lambda|\ge\kappa$.
\item If $\ru(\Sa M)\le\lambda$ then $|S_\lambda|\le 2^{|T|}$.
\end{enumerate}
\end{fact}

In the proof below a ``type" is a one-type.

\begin{proof}
Suppose $\ru(\Sa M)=\lambda$.
If $p_0$ is the restriction of $p\in S_\lambda$ to the empty set then we have $\ru(p)\le\ru(p_0)\le\ru(\Sa M)$, hence $\ru(p_0)=\lambda$, hence $p$ is a non-forking extension of $p_0$.
There are $\le 2^{|T|}$-types over $\varnothing$ and every type over $\varnothing$ has at most $2^{|T|}$ nonforking extensions~\cite[Prop~2.20(iv)]{pillay-book}.
If $\ru(\Sa M)<\lambda$ then $S_\lambda$ is empty.
So we see that $|S_\lambda|\le 2^{|T|}$ when $\ru(\Sa M)\le\lambda$.
Suppose that $\ru(\Sa M)>\lambda$.
Then there is a type $p$ over $\varnothing$ such that $\ru(p)>\lambda$.
By saturation and the definition of U-rank there is a set $A$ of parameters from $\Sa M$ such that $p$ has $\ge\kappa$ extensions over $A$ of U-rank $\ge\lambda$, and each extension $q$ has a global non-forking extension $q^*$, and such $q^*$ is necessarily in $S_\lambda$.
Hence $|S_\lambda|\ge\kappa$.
\end{proof}

Let $\kappa$ range over cardinals $\ge 2$ and let $\tensor[^{\omega>\hspace{.05cm}}]{\kappa}{}$ be set of finite sequences of elements of $\kappa$.
Let $\Sa T_\kappa$ be the unary relational structure on $\tensor[^{\omega\hspace{.05cm}}]{\kappa}{}$ with relations $(U_\sigma : \sigma \in \tensor[^{\omega>\hspace{.05cm}}]{\kappa}{})$ given by
\[
 U_\sigma(\beta) \quad\Longleftrightarrow\quad \sigma \text{  is an initial segment of  } \beta \text{  for all}\quad \sigma \in \tensor[^{\omega>\hspace{.05cm}}]{\kappa}{}, \beta \in \tensor[^{\omega\hspace{.05cm}}]{\kappa}{}.
 \]
In particular $\Sa T_\omega$ is the structure $\Sa U$ described above.

\begin{proposition}
\label{prop:supe}
The following are equivalent:
\begin{enumerate}
\item $T$ is not superstable.
\item $T$ trace defines $\Sa T_\kappa$ for $\kappa=\left(2^{|T|}\right)^+$.
\item $T$ trace defines $\Sa T_\kappa$ for all cardinals $\kappa$.
\item $T$ trace defines any unary relational structure whose unary relations form a well-founded tree of height $\omega$ under inclusion.
\end{enumerate}
So a countable theory is superstable if and only if it does not trace define $\Sa T_{\kappa}$ for $\kappa=\left(2^{{\aleph_0}}\right)^+$.
\end{proposition}

It is easy to see that (3) implies (1).
Note that the language of $\Sa T_\kappa$ has cardinality $\kappa$ and $|S_1(\Sa T_\kappa)| = \kappa^{\aleph_0}$.
Proposition~\ref{prop:unary char} shows that if $T$ trace defines $\Sa T_\kappa$ then  there is $\Sa M\models T$, $n\ge 1$, and $A\subseteq M$ such that $|A|=\kappa$ and $|S_n(\Sa M,A)| \ge \kappa^{\aleph_0}$.
Hence $T$ is not superstable when $T$ trace defines each $\Sa T_\kappa$.
We give another proof below, yielding a better bound on $\kappa$.

\medskip
We let $\sigma \frown \eta$ be the concatenation of finite sequences $\sigma, \eta$.

 
\begin{proof}
We first show that (1) implies (3).
Suppose that $T$ is not superstable and fix a cardinal $\kappa$.
Let $\Sa M\models T$ be $\kappa^+$-saturated.
By \cite[Thm~6.4]{shelah-superstable} there is $n \ge 1$, and a collection $(X_\sigma : \sigma\in\tensor[^{\omega>\hspace{.05cm}}]{\kappa}{})$ of definable subsets of $M^n$ such that
\begin{enumerate}
\item $X_\sigma$ is contained in $X_\rho$ when $\sigma$ extends $\rho$.
\item $X_{\sigma\frown i}$ and $X_{\sigma\frown j}$ are disjoint for all $\sigma\in\kappa^{<\omega}$ and distinct $i,j<\kappa$.
\end{enumerate}
Now follow the proof of Proposition~\ref{prop:unary}.

\medskip
It is clear that (3) implies (2).
We show that (2) implies (1).
Let $\kappa=\left(2^{{|T|}}\right)^+$.
Suppose that $\Sa M\models T$ trace defines $\Sa T_\kappa$ via an injection $\uptau \hookrightarrow \tensor[^{\omega\hspace{.05cm}}]{\kappa}{}\to M^n$.
We may suppose that $\tensor[^{\omega\hspace{.05cm}}]{\kappa}{}$ is a subset of $M^n$ and $\uptau$ is the inclusion.
Suppose towards a contradiction that $\Sa M$ is superstable.
For each $\sigma$ fix $\Sa M$-definable $Y_\sigma\subseteq M^n$ such that $U_\sigma = Y_\sigma \cap\tensor[^{\omega\hspace{.05cm}}]{\kappa}{}$.
We may suppose that the U-rank of each $Y_\sigma$ is minimal among definable sets with this property.
After replacing each $Y_\sigma$ with the intersection of all $Y_\eta$ such that $\eta$ is an initial segment of $\sigma$ we suppose that $Y_{\sigma\frown i}\subseteq Y_\sigma$ for all $\sigma\in\tensor[^{\omega>\hspace{.05cm}}]{\kappa}{}$ and $i <\kappa$.
We produce a contradiction by showing that there is a sequence $(\beta_1,\beta_2,\ldots)\in\tensor[^{\omega\hspace{.05cm}}]{\kappa}{}$ such that $\ru(Y_{\beta_1,\ldots,\beta_k,\beta_{k+1}})<\ru(Y_{\beta_1,\ldots,\beta_k})$ for all $k$.
By induction it is enough to fix $\sigma\in\tensor[^{\omega>\hspace{.05cm}}]{\kappa}{}$ and produce $i<\kappa$ such that $\ru(Y_{\sigma\frown i})<\ru(Y_\sigma)$.

\medskip
Let $\lambda=\ru(Y_\sigma)$.
Suppose towards a contradiction that $\ru(Y_{\sigma\frown i})=\lambda$ for all $i < \kappa$.
Suppose that $i_1,\ldots,i_m,j<\kappa$ are distinct.
Then $U_{\sigma\frown j}$ is contained in $Y_{\sigma\frown j}\setminus [Y_{\sigma\frown i_1}\cup\cdots\cup Y_{\sigma\frown i_m}]$, so by minimality we have
\[
\ru(Y_{\sigma\frown j}\setminus [Y_{\sigma\frown i_1}\cup\cdots\cup Y_{\sigma\frown i_m}])=\ru(Y_{\sigma\frown j})=\lambda
\]
Therefore for all $j<\kappa$ there is a complete type $p_j$ over $\Sa M$ such that:
\begin{enumerate}
\item $p_j$ is concentrated on $Y_{\sigma\frown j}\setminus Y_{\sigma\frown i}$ for all $i\ne j$,
\item $p_j$ is not concentrated on any definable set of U-rank $<\lambda$.
\end{enumerate}
We have constructed a set $\{p_j:j<\kappa\}$ of $\kappa$ distinct types of U-rank $\lambda$ each of which is concentrated on $Y_\sigma$.
This is a contradiction by choice of $\kappa$ and Fact~\ref{fact:U rank}(2).

\medskip
It is clear that (4) implies (3).
Let $\Sa X$ be a unary relational structure whose unary relations form a well-founded tree of height $\omega$ under inclusion.
Then $\Cal B[\Sa X]$ embeds into $\Cal B[\Sa T_\kappa]$ for some $\kappa$.
Hence an application of Lemma~\ref{lem:both ways} shows that (3) implies (4).
\end{proof}

We now prove a preservation result for U-rank.
We let $\eta \cdot \lambda$ be the natural product of ordinals $\eta,\lambda$, see Appendix~\ref{section:ordinals} for background.

\begin{proposition}
\label{prop:U rank}
If there is an $n$-ary trace definition $\Sa M \rightsquigarrow \Sa N$ then $\ru(\Sa M)\le n\cdot\ru(\Sa N)$.
In particular if  $T^*$ is trace definable in $T$ then $\ru(T^*)<\omega\cdot\ru(T)$ and hence $T^*$ has finite U-rank when $T$ has finite U-rank.
\end{proposition}

In particular $\ru(\Sa M) \le \ru(\Sa N)$ when $\Sa M$ trace embeds into $\Sa N$.
Proposition~\ref{prop:U rank} follows from Proposition~\ref{prop:U rank fine} and  a case of the Lascar inequalities: $\ru_{\Sa N}(N^n)\le n\cdot\ru_{\Sa N}(N)$ \cite[19.2]{Poizat}.

\begin{proposition}
\label{prop:U rank fine}
Suppose that $M \subseteq N^n$ and that the inclusion $M \hookrightarrow N^n$ is a trace definition $\Sa M \rightsquigarrow \Sa N$.
Suppose that $X$ is an $\Sa M$-definable subset of $M^m$ and $Y$ is an $\Sa N$-definable subset of $N^{mn}$ such that $X = Y\cap M^m$.
Then $\ruo(X)\le\rum(Y)$.
\end{proposition}

We first prove a lemma.
Let $\mathrm{Def}(\Sa M)$ be the collection of definable subsets of $M$.
Given an ordinal $\lambda$ we let $\app$ be the equivalence relation on $\mathrm{Def}(\Sa M)$ given by declaring $X\app X'$ if and only if $\rum(X\triangle X')<\lambda$.

\begin{lemma}
\label{lem:app}
Suppose that $T$ is superstable and $\lambda$ is an ordinal.
The following are equivalent:
\begin{enumerate}
\item $\ru(T)\le\lambda$.
\item $|\mathrm{Def}(\Sa M)/\!\app\!|\le 2^{2^{|T|}}$ for all $\Sa M\models T$.
\end{enumerate}
\end{lemma}

\begin{proof}
Suppose $\Sa M\models T$ and let $S_\lambda$ be as in Fact~\ref{fact:U rank}. 
Let $X,X'$ range over definable subsets of $M$.
Then $X\app X'$ if and only if $X$ and $X'$ contain the same types of rank $\ge\lambda$.
Hence $|\mathrm{Def}(\Sa M)/\!\app\!|\le 2^{|S_\lambda|}$.
So if $\ru(T) \le \lambda$ then $|\mathrm{Def}(\Sa M)/\!\app\!|\le 2^{2^{|T|}}$ by Fact~\ref{fact:U rank}.
If $p\in S_\lambda$ and $X\app X'$ then $p$ is not in $X\triangle X'$, so $p$ is in $X$ if and only if $p$ is in $X'$.
Hence we also have $|S_\lambda| \le 2^{|\mathrm{Def}(\Sa M)/\!\app\!|}$.
Again apply Fact~\ref{fact:U rank}.
\end{proof}

\begin{proof}[Proof of Proposition~\ref{prop:U rank fine}]
To simplify notation we only treat the case when $X=M$ and $Y=N$.
The general case follows in the same way by replacing $\Sa M$, $\Sa N$ with the structure induced on $X$, $Y$ by $\Sa M$, $\Sa N$, respectively.
Let $T$, $T^*$ be the theory of $\Sa N$, $\Sa M$, respectively.
After possibly passing to elementary extensions and adding constant symbols we may suppose that $|T|=|T^*|$.
If $\Sa N$ is not superstable then $\ru(\Sa N)=\infty$ and the inequality trivially holds.
We suppose that $\Sa N$ is superstable.
Then $\Sa M$ is superstable by Corollary~\ref{cor:superstable}.

\medskip
Let $\lambda$ be an ordinal.
We apply induction on $\lambda$ to show that if $\ru(\Sa N)\le\lambda$ then $\ru(\Sa M)\le\lambda$.
Suppose $\lambda=0$.
Then $\Sa N$ is finite, hence $\Sa M$ is finite, hence $\ru(\Sa M) = 0$.
Suppose $\lambda > 0$.
By Lemma~\ref{lem:app} and Proposition~\ref{prop:trace-theories} it is enough to show that $|\mathrm{Def}(\Sa M)/\!\app\!|\le 2^{2^{|T^*|}}$.
By Lemma~\ref{lem:app} it is enough to show that $|\mathrm{Def}(\Sa M)/\!\app\!|\le|\mathrm{Def}(\Sa N)/\!\app\!|$.
It is enough to fix $\Sa M$-definable $X,X' \subseteq M$ and $\Sa N$-definable $Y,Y'\subseteq N$, suppose that $X=Y\cap M$, $X'=Y'\cap M$ and $Y\app Y'$, and show that $X \app X'$.
We have $\ru_{\Sa N}(Y\triangle Y')<\lambda$ as $Y\app Y'$.
Note that $X\triangle X' \subseteq Y\triangle Y'$, so $\ru_{\Sa M}(X\triangle X')<\lambda$ by induction, hence $X\app X'$.
\end{proof}

We now give preservation results for Morley rank.

\begin{proposition}
\label{prop:morley rank}
If $\Sa N$  trace defines $\Sa M$ then $\mr(\Sa M) < \mr(\Sa N)^\omega$.
If $T$ has finite Morley rank then any theory trace definable in $T$ has finite Morley rank.
\end{proposition}

Proposition~\ref{prop:morley rank} follows from Proposition~\ref{prop:morely rank fine} and the fact that $\mr_{\Sa N}(N^n)\le [\mr(\Sa N)+1]^n$ for all $n \ge 1$, see \cite[Thm~V 7.8]{ShelahCT} or \cite[Exercise~6.4.4]{Tent_Ziegler_2012}.
(Here the power is the usual ordinal power.)

\begin{proposition}
\label{prop:morely rank fine}
Suppose that $M \subseteq N^n$ and that the inclusion $M \hookrightarrow N^n$ is a trace definition $\Sa M\rightsquigarrow \Sa N$.
Suppose that $X$ is an $\Sa M$-definable subset of $M^m$ and $Y$ is an $\Sa N$-definable subset of $N^{mn}$ such that $X = Y \cap M^m$.
Then $\mr_{\Sa M}(X) \le \mr_{\Sa N}(Y)$. 
\end{proposition}

It follows in particular that $\mr(\Sa M) \le \mr(\Sa N)$ when $\Sa M$ trace embeds into $\Sa N$.



\begin{proof}
To simplify notation we only treat the case when $m = n = 1$, the general case follows in the same way.
Let $\lambda$ be an ordinal and suppose that $\mrm(X) \ge \lambda$.
We show that $\mr_{\Sa N}(Y) \ge \lambda$.
We apply induction on $\lambda$.
If $\mrm(X) \ge 0$ then $X$ is nonempty, hence $Y$ is nonempty, hence $\mr_{\Sa N}(Y) \ge 0$.
The case when $\lambda$ is a limit ordinal is clear, we suppose that $\lambda$ is a successor ordinal.
Fix $n$.
By the definition of Morley rank it is enough to produce pairwise disjoint $\Sa N$-definable $Y_1,\ldots,Y_n \subseteq Y$ such that $\mr_{\Sa N}(Y_i) \ge \lambda - 1$ for each $i$.
As $\mrm(X) > \lambda - 1$ there are pairwise disjoint $\Sa M$-definable $X_1,\ldots,X_n \subseteq X$ such that $\mrm(X_i) \ge \lambda - 1$ for each $i$.
For each $i$ fix $\Sa N$-definable $Y_i \subseteq N$ such that $Y_i \cap M = X_i$.
After replacing each $Y_i$ with $Y \cap Y_i$ we suppose that each $Y_i$ is contained in $Y$.
After replacing each $Y_i$ with $Y_i \setminus \bigcup_{j \ne i} Y_j$ we suppose that the $Y_i$ are pairwise disjoint.
By induction we have $\mr_{\Sa N}(Y_i) \ge \mrm(X_i) \ge \lambda - 1$.
\end{proof}

We next give a characterization result for Morley rank.
We first need to prove a lemma.
We say that a unary structure $\Sa X$ is {\bf type-determined} if every $a \in X$ is determined by its type, i.e. if $a \mapsto \tp_{\Sa X}(a)$ gives an injection $X \hookrightarrow S_1(\Sa X)$.

\begin{lemma}\label{lem:type determined}
Let $\Sa X$ be a type-determined unary structure, $\Sa M$ an $|L|^+$-saturated $L$-structure, and $f \colon S_1(\Sa X) \hookrightarrow S_1(\Sa M)$ be a continuous embedding.
Then $\Sa X$ trace embeds into $\Sa M$.
\end{lemma}

\begin{proof}
By saturation every type in $S_1(\Sa M)$ has a realization in $\Sa M$.
Hence there is a map $\uptau \colon X \to M$ such that $\tp_{\Sa M}(\uptau(a))$ realizes $f(\tp_{\Sa X}(a))$ for all $a \in X$.
Now $\uptau$ is injective as $\Sa X$ is type-determined.
The proof of Lemma~\ref{lem:both ways} shows that $\uptau$ is a trace embedding $\Sa X \hookrightarrow \Sa M$.
\end{proof}

Given an ordinal $\lambda$ let $\Sa X_\lambda$ be the unary relational structure with domain $\omega^\lambda$ and unary relations defining $[0,\eta]$ for every $\eta<\omega^\lambda$.
Note that $\Sa X_\lambda$ is type-determined.

\begin{proposition}
\label{prop:unary 2}
Fix an ordinal $1\le\lambda<\omega_1$.
Suppose that $\Sa M$ is an $\aleph_1$-saturated.
Then $\lambda \le \mr(\Sa M)$ if and only if there is a trace embedding $\Sa X_\lambda \hookrightarrow \Sa M$.
Hence  $\Sa X_\lambda$ is trace definable in $\Sa M$ if and only if $\lambda \le \mr_{\Sa M}(M^n)$ for some $n$.
\end{proposition}

Recall that $\CB(X)$ is the Cantor rank of a compact Hausdorff space $X$.

\begin{proof}
The second claim is immediate from the first.
We prove the first claim.
Equip $\omega^\lambda + 1$ with the order topology and note that this topology is homeomorphic to $S_1(\Sa X_\lambda)$.
Then $\mr(\Sa X_\lambda) \ge \mathrm{CR}(\omega^\lambda + 1) = \lambda$.
So the right to left implication follows by Proposition~\ref{prop:morely rank fine}.
We prove the other implication.
Suppose $\mr(\Sa M) \ge \lambda$.
Let $L$ be the language of $\Sa M$.
We only need countably many definable sets to witness $\mr(\Sa M)\ge\lambda$ so some reduct of $\Sa M$ to a countable sublanguage has Morley rank $\ge \lambda$.
We suppose that the language of $\Sa M$ is countable after possibly passing to this sublanguage.
After possibly adding countably many constants to the language we suppose that the witness to $\mr(\Sa M)\ge\lambda$ is zero-definable.
Therefore  $S_1(\Sa M)$ has Cantor rank $\ge\lambda$.
Note that $S_1(\Sa M)$ is second countable.
Hence Fact~\ref{fact:db} shows that $S_1(\Sa X_\lambda) = \omega^\lambda+1$ continuously embeds into $S_1(\Sa M)$.
Apply Lemma~\ref{lem:type determined}.
\end{proof}

Corollary~\ref{cor:unary} is immediate from Proposition~\ref{prop:unary 2}.

\begin{corollary}
\label{cor:unary}
An arbitrary theory $T$ has finite Morley rank if and only if $T$ does not trace define $\Sa X_\omega$.
\end{corollary}

We next classify countable unary theories modulo trace equivalence.
We first prove a lemma about Morley ranks of unary structures.

\begin{lemma}
\label{lem:unary rank}
Suppose $\Sa X$ is a unary structure.
Then $\CB(S_1(\Sa X))\le\mr(\Sa X)\le\CB(S_1(\Sa X))+1$ and $\mr_{\Sa X}(X^n)=n\cdot\mr(\Sa X)$ for all $n\ge 2$.
\end{lemma}

Here $\cdot$ is the Hessenburg product of ordinals as above.

\begin{proof}
After possibly passing to an elementary extension we suppose that $\Sa X$ is highly saturated.
By the definition of Morley rank we have $\CB(S_1(\Sa X))\le\mr(\Sa X)$.
By saturation there is a set $A$ of parameters from $\Sa M$ such that $\mr_{\Sa X}(X^n)=\CB(S_n(\Sa X,A))$ for all $n$.
Now every $A$-definable subset of $X$ is zero-definable modulo a finite subset of $A$.
It follows that $\mr(\Sa X)\le\mathrm{CR}(S_1(\Sa X))+1$.

\medskip
Now suppose that $n\ge 2$.
By Lemma~\ref{lem:db} we have 
\[
n \cdot \mr(\Sa X) = n \cdot \mathrm{CR}(S_1(\Sa X, A)) = \mathrm{CR}(S_1(\Sa X, A)^n).
\]
It suffices to show that $\CB(S_n(\Sa X,A)) = \cra(S_1(\Sa X, A)^n)$.
Now $\CB(S_n(\Sa X,A)) \ge \cra(S_1(\Sa X, A)^n)$ as there is an obvious continuous surjection $S_n(\Sa X, A) \to S_1(\Sa X, A)^n$.
Let $p_1, \ldots, p_d$ be an enumeration of the $n$-types in the theory of equality.
For each $i = 1, \ldots, d$, let $Y_i$ be the set of realizations of $p_i$ in $X^n$, and let $Y^*_i$ be the set of $p \in S_n(\Sa X, A)$ concentrated on $Y_i$.
It is enough to fix $i$ and show that $\cra(Y^*_i) \le \cra(S_1(\Sa X, A)^n)$.
Every $\Sa X$-definable subset of $Y_i$ is of the form $Y_i \cap [E_1 \times \cdots \times E_n]$ for $A$-definable $E_1, \ldots, E_n \subseteq X$ as $\Sa X$ is unary.
Hence the natural map $Y^*_i \to S(\Sa X, A)^n$ is a continuous embedding, which gives the desired inequality.
\end{proof}

We now classify unary structures in countable languages modulo trace equivalence.
We let $\Val$ be the archimedean valuation on Conway's ordered field of surreal numbers $\no$, so if $\xi$ and $\zeta$ are ordinals then $\Val(\xi)\le\Val(\zeta)$ if and only if $\xi\le n\cdot\zeta$ for some $n \ge 1$.
The additive group of $\no$ is canonically identified with the value group of $\Val$ via the map $\gamma \mapsto \Val(\omega^\gamma)$~\cite[Chapter~5.B]{gonshor}.
If $\xi$ is an ordinal with Cantor normal form $\sum_{\lambda\in\on} n_\lambda \cdot \omega^\lambda$ then $\Val(\xi)$ is, under the canonical identification, the maximal ordinal $\lambda$ such that $n_\lambda\ne 0$.
Set $\Val(\infty)=\infty$.

\begin{proposition}
\label{prop:unary c}
Suppose that $\Sa X$ and $\Sa X^*$ are unary structures in countable languages and declare $\xi=\mathrm{RM}(\Sa X)$ and $\xi^*=\mathrm{RM}(\Sa X^*)$.
Then $\Th(\Sa X)$ trace defines $\Th(\Sa X^*)$ if and only if $\Val(\xi^*)\le\Val(\xi)$ and $\Sa X$ is trace equivalent to $\Sa X^*$ if and only if $\Val(\xi)=\Val(\xi^*)$.
Hence exactly one of the following holds:
\begin{enumerate}[leftmargin=*]
\item $\Sa X$ is trace equivalent to  $\Sa C$.
\item $\Sa X$ is trace equivalent to $\Sa X_{\omega^\lambda}$ for a unique ordinal $\lambda$, and this $\lambda$ is the leading exponent in the Cantor normal form of $\mr(\Sa X)$.
\end{enumerate}
\end{proposition}

\begin{proof}
The second claim follows easily from the first claim.
For the first claim it is enough to show that $\Val(\xi^*)\le\Val(\xi)$ if and only if $\Th(\Sa X)$ trace defines $\Th(\Sa X^*)$.
The right to left direction follows by Proposition~\ref{prop:morely rank fine} and Lemma~\ref{lem:unary rank}.
Suppose that $\Val(\xi^*) \le \Val(\xi)$.
If $\Val(\xi) = \infty$ then $\Sa X$ is not totally transcendental, hence $\Th(\Sa X)$ trace defines $\Sa X^*$ by Proposition~\ref{prop:unary}.
Suppose that $\Val(\xi) < \infty$.
Fix $n$ such that $\xi^*< n\cdot\xi$.
By Lemma~\ref{lem:unary rank} we have $\xi^*<\mr_{\Sa X}(X^n)$.
Let $\lambda=\xi^*+1$, so $\xi^*<\lambda\le\mr_{\Sa X}(X^n)$.
By Proposition~\ref{prop:unary 2} $\Th(\Sa X)$ trace defines $\Sa X_\lambda$.
By Fact~\ref{fact:db} there is a continuous injection $S_1(\Sa X^*) \hookrightarrow S_1(\Sa X_\lambda)$.
Hence $\Th(\Sa X_{\lambda})$ trace defines $\Sa X^*$ by Lemma~\ref{lem:both ways}.
\end{proof}


\begin{corollary}
\label{cor:unmin}
The following are equivalent for any unary structure $\Sa X$ in a countable language.
\begin{enumerate}
\item $\Sa X$ has finite Morley rank.
\item $\Sa X$ is trace definable in the theory of an infinite set with equality.
\end{enumerate}
\end{corollary}

\begin{proof}
Suppose that $\Sa X$ has finite Morley rank.
Then $\Sa X$ is trace equivalent to either $\Sa X_0$ or $\Sa X_1$.
In the first case $\Sa X$ is finite.
In the second case $\Sa X$ is trace equivalent to an infinite set equipped with equality.
Proposition~\ref{prop:morley rank} gives the other direction.
\end{proof}

We have shown that trace equivalence classes of unary structures in countable languages form a linear order of order type $\omega_1+1$.
We show in Corollary~\ref{cor:non linear} below that this does not extend to languages of cardinality $\aleph_1$.

\section{Dp-rank and op-dimension}
We assume familiarity with dp-rank.
See \cite[Chapter~4]{Simon-Book} for the definition and basic facts.
We let $\dprk(T)$ be the dp-rank of a theory $T$, let $\dprk_{\Sa M}(X)$ be the dp-rank of a definable set $X$ in a structure $\Sa M$, and let $\dprk(\Sa M) = \dprk_{\Sa M}(M)$.

\begin{proposition}
\label{thm:dp-rank}
Suppose that $T^*$ is trace definable in $T$.
If $T$ has finite dp-rank then $T^*$ has finite dp-rank.
If $T$ has infinite dp-rank then $\dprk(T^*) \le \dprk(T)$.
In particular if $T$ is strongly dependent then $T^*$ is strongly dependent.
\end{proposition}

The first two claims of Proposition~\ref{thm:dp-rank} are immediate from Proposition~\ref{prop:dp-rank}.
The last claim follows from the previous as $T$ is strongly dependent if and only if $\dprk(T) < \aleph_0$.
(This is not equivalent to the dp-rank being finite, see \cite[\S 4.2, pg. 54]{Simon-Book}.)

\begin{proposition}
\label{prop:dp-rank}
Suppose that there is an $n$-ary trace definition $\Sa M \rightsquigarrow \Sa N$.
If $\Sa N$ has finite dp-rank then $\dprk(\Sa M) \le n \dprk(\Sa N)$ and if $\Sa N$ has infinite dp-rank then $\dprk(\Sa M) \le \dprk(\Sa N)$.
\end{proposition}

In particular $\dprk(\Sa M) \le \dprk(\Sa N)$ when $\Sa M$ trace embeds into $\Sa N$.
Proposition~\ref{prop:dp-rank} follows from Proposition~\ref{prop:dp-rank-0} and subadditivity of dp-rank.

\begin{proposition}
\label{prop:dp-rank-0}
Suppose that $M \subseteq N^n$ and that the inclusion $M \hookrightarrow N^n$ is a trace definition $\Sa M \rightsquigarrow \Sa N$.
Suppose that $X$ is an $\Sa M$-definable subset of $M^m$ and $Y$ is an $\Sa N$-definable subset of $N^{mn}$ such that $X = Y \cap M^m$.
Then $\dprk_{\Sa M}(X) \le \dprk_{\Sa N}(Y)$. 
\end{proposition}

\begin{proof}
We only treat the case when $m = n = 1$, the general case follows in the same way.
Suppose that $\lambda$ is a cardinal and $\dprk_{\Sa M} X \ge \lambda$.
By Proposition~\ref{prop:trace-theories} we may suppose that $\Sa N \models T$ and $\Sa M \models T^*$ are both highly saturated.
Fix a sequence  of parameter free $L^*$-formulas $(\varphi_\alpha(x_\alpha , y) : \alpha < \lambda )$ and an array of tuples $(a_{\alpha,i} \in M^{|x_\alpha|} : \alpha < \lambda , i < \omega )$ such that for any function $f \colon \lambda \to \N$ there is $b \in X$ such that
$$ \Sa M \models \varphi_\alpha(a_{\alpha,k},b) \quad \Longleftrightarrow \quad f(\alpha) = k \quad \text{for all   } \alpha<\lambda,i\in\N. $$
Fix an $L(N)$-formula $\theta_\alpha(z_\alpha,w)$ for each $\alpha < \lambda$ with $|z_\alpha| = |x_\alpha|$ and $|w| = 1$, such that for any $a \in M^{|x_\alpha|}$, $b \in M$ we have $\Sa M \models \varphi_\alpha(a,b)$ if and only if $\Sa N \models \theta_\alpha(a,b)$.
Consideration of $(\theta_\alpha(z_\alpha,w) : \alpha < \lambda)$ and $( a_{\alpha,i} : \alpha < \lambda, i < \omega )$ shows that $\dprk_{\Sa N}(Y) \ge \lambda$.
Hence we have $\dprk_{\Sa N}(Y) \ge \dprk(\Sa M)$.
\end{proof}


We now characterize dp-rank in terms of unary structures.
Given sets $A,B$ we let $\tensor[^{A\hspace{.05cm}}]{B}{}$ be the set of functions $A \to B$.
Given cardinals $\lambda$ and $\kappa$ we let $\Sa B^\lambda_\kappa$ be the unary relational structure with domain $\tensor[^{\kappa\hspace{.05cm}}]{\lambda}{}$ and a unary relation $P_{ij}$ for each $i<\kappa,j<\lambda$ given by declaring $P_{ij}(\sigma)$ if and only if $\sigma(i)=j$.

\begin{proposition}
\label{prop:dpdp}
Suppose that $\kappa$ is a cardinal and $\Sa M$ is $\lambda^+$-saturated for  $\lambda = \kappa + |T|$.
Then the following are equivalent.
\begin{enumerate}
\item $\kappa\le\dprk(\Sa M)$.
\item There is a trace embedding $\Sa B^\lambda_\kappa \hookrightarrow \Sa M$.
\item There is a trace embedding $\Sa B^\xi_\lambda \hookrightarrow \Sa M$ for $\xi=|T|^+$.
\end{enumerate}
\end{proposition}


\begin{proof}
Suppose that $\dprk(\Sa M)\ge\kappa$.
By saturation there are formulas $(\varphi_i(x,y_i) : i<\kappa)$ and tuples of parameters $(\beta_{ij} : i<\kappa,j<\lambda)$ and $(\gamma_\sigma : \sigma\in\tensor[^{\kappa\hspace{.05cm}}]{\lambda}{})$ from $\Sa M$ such that
$$ \Sa M\models \varphi_i(\gamma_\sigma,\beta_{ij})\quad\Longleftrightarrow\quad \sigma(i)=j \quad\text{  for all } i<\kappa,j<\lambda,\sigma\in\tensor[^{\kappa\hspace{.05cm}}]{\lambda}{}.$$
Let $\uptau\colon\tensor[^{\kappa\hspace{.05cm}}]{\lambda}{} \to M^{|x|}$ be given by $\uptau(\sigma)=\gamma_\sigma$.
Then $\uptau$ is an injection and we have 
$$ \Sa B^\lambda_\kappa\models P_{ij}(\sigma)\quad\Longleftrightarrow\quad \Sa M\models\varphi_i(\uptau(\sigma),\beta_{ij}) \quad\text{  for all }  i<\kappa,j<\lambda,\sigma\in\tensor[^{\kappa\hspace{.05cm}}]{\lambda}{}$$
Hence $\uptau$ is a trace embedding $\Sa B^\lambda_\kappa \hookrightarrow \Sa M$.
Hence (2) holds.

\medskip
Note that (2) implies (3).
We show that (3) implies (1).
Let $\xi=|T|^+$ and suppose that $\uptau$ is a trace embedding $\Sa B^\xi_\kappa \hookrightarrow \Sa M$.
We show $\dprk_{\Sa M}(M)\ge\kappa$.
For all $i<\kappa$, $j<\xi$ let $Y_{ij}\subseteq M$ be $\Sa M$-definable such that $\uptau^{-1}(Y_{ij})=P_{ij}$.
Fix $i<\kappa$.
As $\xi>|T|$ there is a parameter-free formula $\varphi_i$ such that infinitely many $Y_{ij}$ may be defined using instances of $\varphi_i$.
Hence there is a sequence $(\beta_{ij} : j <\omega)$ of parameters and an injection $f_i\colon\omega\to \xi$ such that $\varphi_i(x,\beta_{ij})$ defines $Y_{i f_i(j)}$ for all $j<\omega$.
Now select such $\varphi_i$, $(\beta_{ij}:j<\omega)$, and $f_i\colon\omega\to\xi$ for all $i<\kappa$.
Let $g\colon\tensor[^{\kappa\hspace{.05cm}}]{\omega}{} \to \tensor[^{\kappa\hspace{.05cm}}]{\xi}{}$ be given by declaring $g(\sigma)=\eta$ when $\eta(i)=f_i(\sigma(i))$ for all $i<\kappa$.
Let $\gamma_\sigma$ be $\uptau(g(\sigma))$ for all $\sigma\in\tensor[^{\kappa\hspace{.05cm}}]{\omega}{}$.
Now observe that
$$ \Sa M\models \varphi_i(\gamma_\sigma,\beta_{ij})\quad\Longleftrightarrow\quad \sigma(i)=j \quad\text{  for all } i<\kappa,j<\omega,\sigma\in\tensor[^{\kappa\hspace{.05cm}}]{\omega}{}.$$
Hence $\dprk_{\Sa M}(M)\ge\kappa$.
\end{proof}

\begin{proposition}
\label{prop:undep}
Let $\kappa$ be an infinite cardinal.
The following are equivalent:
\begin{enumerate}
\item The dp-rank of $T$ is at least $\kappa$.
\item $T$ trace defines $\Sa B^\lambda_\kappa$ for $\lambda=|T|^+$.
\item $T$ trace defines $\Sa B^\lambda_\kappa$ for every cardinal $\lambda$.
\end{enumerate}
\end{proposition}

\begin{proof}
Proposition~\ref{prop:dpdp} shows that (2) and (3) are equivalent and either holds if and only if $\dprk_{\Sa M}(M^n) \ge \kappa$ for some $n$.
An application of subadditivity for dp-rank shows that $\dprk_{\Sa M}(M^n)\ge\kappa$ implies $\dprk_{\Sa M}(M)\ge\kappa$.
\end{proof}

Corollary~\ref{cor:undep} follows  as $T$ is strongly dependent if and only if $\dprk(T)<\aleph_0$.

\begin{corollary}
\label{cor:undep}
The following are equivalent for any theory $T$.
\begin{enumerate}
\item $T$ is not strongly dependent.
\item $T$ trace defines $\Sa B^\lambda_\omega$ for $\lambda=|T|^+$.
\item $T$ trace defines $\Sa B^\lambda_\omega$ for every cardinal $\lambda$.
\end{enumerate}
In particular a countable theory is not strongly dependent if and only if it trace defines $\Sa B^{\omega_1}_\omega$.
\end{corollary}

Let $\Sa O$ be the unary relational structure with domain $\omega_1$ and relations defining $[0, \eta]$ for every $\eta < \omega_1$.
We have shown that a countable theory $T$ trace defines $\Sa B^{\omega_1}_\omega$, $\Sa O$ if and only if $T$ is not strongly dependent, unstable, respectively. 
It follows that $\Sa B^{\omega_1}_\omega$ and $\Sa O$ are incomparable under trace definability as the classes of countable strongly dependent theories and countable stable theories are incomparable under containment.
Corollary~\ref{cor:non linear} follows.

\begin{corollary}\label{cor:non linear}
Trace equivalence classes of unary structures in languages of cardinality $\aleph_1$ do not form a linear order under trace definability.
\end{corollary}

We next characterize finiteness of dp-rank.
This requires joins of unary structures.
Let $I$ be an index set, $(\Sa X_i)_{i \in I}$ be a family of unary structures, and $L_i$ be the language of $\Sa X_i$ for each $i \in I$.
We construct the join $\Sa X_\sqcup$ of the $\Sa X_i$; this is a unary structure.
We suppose that each $\Sa X_i$ is unary relational, so in general the join is well defined up to interdefinability.
The domain $X_\sqcup$ of $\Sa X_\sqcup$ is the disjoint union of the $X_i$.
The language of $\Sa X_\sqcup$ consists of equality, the disjoint union of all unary relations in the $L_i$, and a new unary relation $P_i$ for each $i \in I$.
Each $P_i$ defines $X_i$ as a subset of $X_\sqcup$.
For any unary $R \in L_i$ and $a \in X_\sqcup$ we declare $\Sa X_\sqcup \models R(a)$ if and only if $a \in X_i$ and $\Sa X_i \models R(a)$.

\begin{lemma}\label{lem:disjoint union}
Let $(\Sa X_i)_{i \in I}$ and $\Sa X_\sqcup$ be as above and let $\Sa M$ be an arbitrary structure.
Then $\Sa M$ trace defines $\Sa X_\sqcup$ if and only if $|I| \le |M|$ and there is $n\ge 1$ such that there is an $n$-ary trace definition $\Sa X_i \rightsquigarrow \Sa M$ for each $i \in I$. 
\end{lemma}

\begin{proof}
The left to right implication is easy and left to the reader.
Suppose that the right hand side holds for $n \ge 1$.
For each $i \in I$ fix an $n$-ary trace definition $\uptau_i \colon \Sa X_i \rightsquigarrow \Sa M$.
Furthermore fix an injection $\sigma \colon I \hookrightarrow M$.
An application of Lemma~\ref{lem:un-trace} shows that the map $X_\sqcup \hookrightarrow M \times M^n$ given by $a \mapsto (\sigma(i), \uptau_i(a))$ for any $i \in I$ and $a \in X_i$ is a trace definition $\Sa X_\sqcup \rightsquigarrow \Sa M$.
\end{proof} 
\begin{corollary}\label{cor:findp}
The following are equivalent for any theory $T$.
\begin{enumerate}
\item $T$ has infinite dp-rank.
\item $T$ trace defines $\bigvee_{n\ge 1}\Sa B^\lambda_n$ for $\lambda=|T|^+$.
\item $T$ trace defines $\bigvee_{n\ge 1}\Sa B^\lambda_n$  for every cardinal $\lambda$.
\end{enumerate}
\end{corollary}

\begin{proof}
Set $\Sa B^\lambda=\bigvee_{n\ge 1}\Sa B^\lambda_n$.
It suffices to show that (1) implies (3) and (2) implies (1).
Suppose that (1) holds, let $\lambda$ be an infinite cardinal, and suppose that $\Sa M\models T$ is $\lambda^+$-saturated.
We have $n\le \dprk(\Sa M)$ for all $n$ so by Proposition~\ref{prop:dpdp} $\Sa B^\lambda_n$ trace trace embeds into $\Sa M$ for every $n\ge 1$.
Hence $\Sa M$ trace defines $\Sa B^\lambda$ by Lemma~\ref{lem:disjoint union}.

\medskip
Now suppose that (2) holds and set $\lambda=|T|^+$.
By Lemma~\ref{lem:disjoint union} there is $m$ such that $\Sa M$ trace defines every $\Sa B^\lambda_n$ via an injection into $M^m$.
By Proposition~\ref{prop:dpdp} we have $\dprk_{\Sa M}(M^m)\ge n$ for all $n$.
Subadditivity of dp-rank shows that $\dprk(\Sa M)$ is infinite.
\end{proof}

\subsection{Op-dimension}\label{section:op}
We now consider op-dimension.
This was introduced by Guingona and Hill~\cite{gh-op}.
Let $\opd_{\Sa M}(X)$ be the op-dimension of an $\Sa M$-definable set $X$, $\opd(\Sa M) = \opd_{\Sa M}(M)$, and $\opd(T) = \opd(\Sa M)$ for any $\Sa M \models T$.
A $k$-order is a set equipped with $k$ linear orders.
Finite $k$-orders form a Ramsey class~\cite[Cor.~1.4]{Bodirsky2014}.
We let $\Sa O_k$ be the generic $k$-order.

\begin{fact}
\label{fact:opd}
Suppose that $X$ and $Y$ are definable sets over $\Sa M\models T$.
\begin{enumerate}[leftmargin=*]
\item\label{o2} If $\opd(T)$ is finite then $T$ is $\nip$.
\item\label{o3} $\opd(X\times Y)\le\opd(X)+\opd(Y)$.
\item\label{o6} If $T$ is $\nip$ and $\monster$ is a monster model of $T$ then $T$ has op-dimension $\ge k$ if and only if there is an uncollapsed indiscernible picture $O \to \monsterset$ of $\Sa O_k$ in $\monster$ over some small set of parameters.
\end{enumerate}
\end{fact}

See \cite[Section 3.1]{gcs} for (\ref{o2}),  \cite[Thm~2.2]{gh-op} for (\ref{o3}), and \cite[Thm~3.4]{gcs} for (3).

\begin{proposition}
\label{prop:op}
If there is an $n$-ary trace definition $\Sa M \rightsquigarrow \Sa N$ then $\opd(\Sa M) \le n \opd(\Sa N)$.
If $\kappa$ is an infinite cardinal then the collection of theories of op-dimension $<\kappa$ is closed under trace definability, hence theories of finite op-dimension are closed under trace definability.
\end{proposition}

\begin{proof}
The second claim follows from the first.
We prove the first claim.
By subadditivity of op-dimension it is enough to suppose that there is an $n$-ary trace definition $\Sa M \rightsquigarrow \Sa N$ and show that $\opd(\Sa M) \le \opd_{\Sa N}(N^n)$.
After replacing $\Sa N$ with the induced structure on $N^n$ is it enough to suppose that there is a trace embedding $\Sa M \hookrightarrow \Sa N$ and show that $\opd(\Sa M) \le \opd(\Sa N)$.
This follows immediately from the definition of op-dimension in terms of IRD-arrays \cite[Def.~1.22]{gh-op} and is left to the reader.
\end{proof}


We give a second characterization of op-dimension.
Let $(I_1;\triangleleft_1),\ldots,(I_n;\triangleleft_n)$ range over linear orders.
We let $\Sa I_\times$ be the unary relational structure with domain $I_\times=I_1\times\cdots\times I_n$ and a unary relation defining the set of $(\alpha_1,\ldots,\alpha_n)\in I_\times$ such that $\alpha_j \triangleleft_j \beta$ for each $j=1,\ldots,n$ and $\beta\in I_j$.
If $(I_1;\triangleleft_1) = \ldots = (I_n ; \triangleleft) = \Sa I$ then we write $\Sa I^n$ for $\Sa I_\times$.

\begin{proposition}
\label{prop;oop}
The following are equivalent for any theory $T$ and $n\ge 1$.
\begin{enumerate}[leftmargin=*]
\item $T$ has op-dimension $\ge n$.
\item The generic $n$-order trace embeds into a model of $T$.
\item $\Sa I_\times$ trace embeds into a model of $T$ for any linear orders $(I_1;\triangleleft_1),\ldots,(I_n;\triangleleft_n)$.
\item $\Sa I^n$ trace embeds into a model of $T$ for some linear order $\Sa I = (I;\triangleleft)$ with $|I|>\beth_{n-1}(|T|)$.
\end{enumerate}
\end{proposition}

We first recall the \Erdos-Rado theorem: If $\kappa$ is an infinite cardinal, $n\ge 2$, $Y$ is a set of cardinality $\kappa$, $X$ is a set of cardinality exceeding $\beth_{n-1}(\kappa)$,  $h$ is a function $X^n\to Y$, and $<$ is a linear order on $X$, then there is a subset $Z\subseteq X$ such that $\kappa < |Z|$ and $h$ is constant on the set of $(a_1,\ldots,a_n) \in Z^n$ such that $a_1 < \cdots < a_n$. 

\begin{proof}
The case $n=1$ follows from Proposition~\ref{prop:stable-0} so we may suppose that $n\ge 2$.
Equivalence of (1) and (2) follows by  Proposition~\ref{prop:picture} and Fact~\ref{fact:opd}(\ref{o6}).
It is clear that (3) implies (4).
We show that (1) implies (3).
Suppose that (1) holds.
Let $\Sa M\models T$ be $\kappa^+$-saturated for $\kappa = |I_1| + \cdots + |I_n|$.
By an easy saturation argument there is an IRD-array $\varphi_1(x;y_1),\ldots,\varphi_n(x;y_n)$, $(b^i_j\in M^{|y_i|}: i\in\{1,\ldots,n\}, j\in I_i\})$ in $\Sa M$.
By the definition of IRD-arrays there is a function $\uptau\colon I_\times \to M$ such that we have
\[
\Sa M\models\varphi_i(\uptau(a_1,\ldots,a_n);b^i_j) \quad\Longleftrightarrow \quad j\triangleleft_i a_i\quad\text{for $i\in\{1,\ldots,n\}$, $a_1\in I_1,\ldots,a_n\in I_n, j\in I_i$.}
\]
It follows by Lemma~\ref{lem:un-trace} that $\uptau$ gives a trace embedding $\Sa I_\times \hookrightarrow \Sa M$.

\medskip
We show that (4) implies (1).
Suppose that $\Sa I$ and $\Sa M\models T$ are as in (3).
We may suppose that $I^n\subseteq M$ and that the  inclusion $\Sa I^n \hookrightarrow \Sa M$ is a trace embedding.
Let $$X_i^\gamma=\{(a_1,\ldots,a_n)\in I^n : a_i\triangleleft \gamma\} \quad\text{for every}\quad \gamma \in I,i\in\{1,\ldots,n\}.$$
Now for each $i=1,\ldots,n$ and $\alpha\in I$ fix a parameter-free formula $\varphi^i_\alpha(x;y_\alpha)$ and tuple $\gamma_\alpha\in M^{|y_\alpha|}$ such that $|x|=1$ and $I^n\cap \varphi_\alpha(M;\gamma_\alpha) = X^{\gamma_\alpha}_i$.
By the \Erdos-Rado theorem there is an infinite $J\subseteq I$ and parameter-free formulas $\vartheta_1,\ldots,\vartheta_n$ such that $$(\varphi^1_{\alpha_1},\ldots,\varphi^n_{\alpha_n})=(\vartheta_1,\ldots,\vartheta_n)\quad \text{for all $\alpha_1<\cdots<\alpha_n$ in $J$}.$$
Then for any $b=(b_1,\ldots,b_n)\in J^n$ and $\alpha\in J$ we have $\Sa M\models \vartheta_i(b;\gamma_\alpha)$ if and only if $b_i\triangleleft \alpha$.
As $J$ is an infinite chain it contains either an infinite increasing sequence or an infinite decreasing sequence.
We treat the first case, the second case follows by a similar argument.
Let $(b_i)_{i\in\N}$ be a strictly increasing sequence of elements of $J$.
Set $b_{i,j}=b_i$ for all $i\in\N$ and $j\in\{1,\ldots,n\}$.
Finally observe that $(b_{i,j}:i\in\N,j\in\{1,\ldots,n\})$ and $\vartheta_1,\ldots,\vartheta_n$ form an IRD-array of depth $n$ and length $\aleph_0$.
Hence $T$ has op-dimension at least $n$.
\end{proof}

\begin{corollary}
\label{cor:op}
The following are equivalent for any theory $T$.
\begin{enumerate}
[leftmargin=*]
\item $T$ has infinite op-dimension.
\item $T$ trace defines $\bigvee_{n\ge 1}\Sa I^n$ for every linear order $\Sa I=(I;<)$.
\item $T$ trace defines $\bigvee_{n\ge 1}\Sa I^n$ for some linear order $\Sa I=(I;<)$ such that $|I|\ge\beth_{\omega}(|T|)$.
\end{enumerate}
\end{corollary}

Corollary~\ref{cor:op} follows from Proposition~\ref{prop;oop} and Lemma~\ref{lem:disjoint union}.
We leave the proof to the reader as it is similar to that of Corollary~\ref{cor:findp}.

\section{Trace maximality}\label{section:tm}
We have seen that certain classification-theoretic properties like $\mathrm{NSOP}$ are not preserved under trace definability.
In this section we show that these properties in fact have no trace-theoretic consequences.
We say that a theory is \textbf{trace maximal} when it trace defines every structure and that a structure is trace maximal when its theory is.

\begin{proposition}\label{prop:psf}
Any pseudofinite field is trace maximal.
More generally any pseudo real closed field which is not separably or real closed is trace maximal.
\end{proposition}

Proposition~\ref{prop:psf} follows from Proposition~\ref{prop:pseudofinite} below.
We first prove a general fact about trace maximality.

\begin{lemma}
\label{lem:max}
Then the following are equivalent for any theory $T$.
\begin{enumerate}
\item $T$ is trace maximal.
\item There is $\Sa M \models T$ and infinite $A \subseteq M^n$ for some $n$ such that every subset of every $A^m$ is a trace of an $\Sa M$-definable set.
\item There is $\Sa M \models T$ and an infinite $A \subseteq M^n$ such that for every $m$-hypergraph $E$ on $A$ there is an $\Sa M$-definable $X \subseteq M^{mn}$ such that we have $E(a_1,\ldots, a_m)$ if and only if $ (a_1,\ldots,a_m) \in X$ for all $a_1,\ldots,a_m \in A$.

\end{enumerate}
\end{lemma}

Here (3) is useful when dealing with commutative algebraic structures.

\begin{proof}
Suppose that $T$ is trace maximal.
Let $A$ be an infinite set and $\Sa A$ be a structure on $A$ which defines every subset of every $A^m$.
Then $\Th(\Sa A)$ is trace definable in $T$.
Hence (1) implies (2).
We show that (2) implies (1).
It is enough to show that $\Sa A$ is trace maximal.
It suffices to show that if $\Sa B$ is a highly saturated elementary extension of $\Sa A$ and $D$ is a small subset of $B$ then every subset of every $D^m$ is a trace of a $\Sa B$-definable set.
By compactness it suffices to fix $m \ge 1$ and produce a subset $X$ of $A^{m} \times A$ such that for every pair $E, F$ of disjoint finite subsets of $A^m$ there is $b \in A$ such that $E$ is contained in $X_b := \{a \in A^m : (a, b) \in X \}$ and $F$ is disjoint from $X_b$.
This follows as the collection of pairs $E, F$ has cardinality $|A|$.

\medskip
It is clear that (2) implies (3).
We show that (3) implies (2).
Suppose (3).
For simplicity we only treat the case $n = 1$, the general case follows in the same way.
Let $E$ be a graph on $A$ and $(a_i)_{i < \omega }$, $(b_j)_{j < \omega }$ be sequences of distinct elements of $A$ such that for all $i,j$ we have $E(a_i,b_j)$ if and only if $i < j$.
Let $L$ be the language of $\Sa M$ and let $\delta(x,y)$ be an $L(M)$-formula such that for all $a,a^* \in A$ we have $E(a,a^*)$ if and only if $\Sa M \models \delta(a,a^*)$.
For each $i$ let $c_i = (a_i,b_i)$ and let $\phi(x_1,y_1,x_2,y_2)$ be $\delta(x_1,y_2)$.
Then for any $i,j$ we have
\begin{align*}
\Sa M \models \phi(c_i,c_j) &\Longleftrightarrow \Sa M \models \phi(a_i,b_i,a_j,b_j)\\
&\Longleftrightarrow \Sa M \models \delta(a_i,b_j) \\
&\Longleftrightarrow i < j.
\end{align*}
We now show that for any $X \subseteq \omega^m$ there is an $\Sa M$-definable $Y \subseteq M^{2m}$ such that for all $(i_1,\ldots,i_m) \in \omega^m$ we have $(i_1,\ldots,i_m) \in X$ if and only if $(c_{i_1},\ldots,c_{i_m}) \in Y$.
Trace maximality of $\Sa M$ follows by Lemma~\ref{lem:max}(3).
We apply induction on $m \ge 1$.

\medskip
Suppose $m = 1$ and $X \subseteq \omega$.
Let $F$ be a graph on $A$ and $d \in A$ be such that we have $F(a_i,d)$ if and only if $i \in X$ for all $i < \omega$.
Let $\theta(x,y)$ be an $L(M)$-formula such that we have $F(a,a^*)$ if and only if $\Sa M \models \theta(a,a^*)$ for all $a,a^* \in A$.
Let $Y$ be the set of $(a,b) \in M^2$ such that $\Sa M \models \theta(a,d)$.
Then for any $i < \omega$ we have $i \in X$ if and only if $c_i \in Y$.

\medskip
We now suppose that $k \ge 2$ and $X \subseteq \omega^m$.
We let $\omega$ denote the structure $(\omega; <)$.
We let $\qftp_\omega(\imag_1,\ldots,\imag_m)$ be the quantifier free type (equivalently: order type) of $(\imag_1,\ldots,\imag_m) \in \omega^m$ and let $S_m(\omega)$ be the set of quantifier free $m$-types.
Note that $S_m(\omega)$ is finite.
For each $p \in S_m(\omega)$ we fix an $\Sa M$-definable $Y_p \subseteq M^{2m}$ such that $$\qftp_{\omega}(\imag_1,\ldots,\imag_m) = p\quad \Longleftrightarrow\quad (c_{\imag_1},\ldots,c_{\imag_m}) \in Y_p\quad \text{for all}\quad  (\imag_1,\ldots,\imag_m) \in \omega^m.$$
We show that for every $p \in S_m(\omega)$ there is an $\Sa M$-definable subset $X_p$ of $M^{2m}$ such that for any realization $(\imag_1,\ldots,\imag_m) \in \omega^m$ of $p$ we have $(c_{\imag_1},\ldots,c_{\imag_m}) \in X$ if and only if $(c_{\imag_1},\ldots,c_{\imag_m}) \in X_p$.
Then for any $(\imag_1,\ldots,\imag_m) \in \omega^m$ we have 
$$(c_{\imag_1},\ldots,c_{\imag_m}) \in X \quad\Longleftrightarrow\quad (c_{\imag_1},\ldots,c_{\imag_m}) \in \! \bigcup_{p \in S_m(\omega)} (Y_p \cap X_p).$$

We fix $p(x_1,\ldots,x_m) \in S_m(\omega)$ and produce $X_p$.
We first treat the case when $p \models (x_i = x_j)$ for some $i \ne j$.
To simplify notation we suppose that $p \models (x_1 = x_2)$.
Let $X^*$ be the set of $(\imag_1,\ldots,\imag_{m - 1}) \in \omega^{m - 1}$ such that $(\imag_1,\imag_1,\imag_2,\ldots,\imag_{m - 1}) \in X$.
By induction there is an $\Sa M$-definable $Y^* \subseteq M^{2(m - 1)}$ such that for any $(\imag_1,\ldots,\imag_{m - 1}) \in \omega^{m - 1}$  we have $(\imag_1,\ldots,\imag_{m - 1}) \in X^*$ if and only if $(c_{\imag_1},\ldots,c_{\imag_{m- 1}}) \in Y^*$.
Let $X_p$ be the set of $(d_1,\ldots,d_m) \in M^{2m}$ such that $d_1 = d_2$ and $(d_2,\ldots,d_m) \in Y^*$, note that $X_p$ is definable in $\Sa M$.

\medskip
Now suppose that  $p \models (x_i \ne x_j)$ when $i \ne j$.
For any distinct $\imag_1,\ldots,\imag_m < \omega$ there is a unique permutation $\sigma_p$ of $\{1,\ldots,m\}$ with $\qftp_{\omega}(\imag_{\sigma_p(1)},\ldots,\imag_{\sigma_p(m)}) = p$.
Let $H^*$ be the $m$-hypergraph on $\{c_\imag : \imag < \omega \}$ where we have $H^*(c_{\imag_1},\ldots,c_{\imag_m})$ if and only if $\imag_1,\ldots,\imag_m$ are distinct and $(\imag_{\sigma_p(1)},\ldots,\imag_{\sigma_p(m)}) \in X$.
Let $H$ be the $m$-hypergraph on $\{a_\imag : \imag < \omega\}$ where we have $H(a_{\imag_1},\ldots,a_{\imag_m})$ if and only if $H^*(c_{\imag_1},\ldots,c_{\imag_m})$ for all $(\imag_1,\ldots,\imag_m) \in \omega^m$.
By assumption there is an $L(M)$-formula $\varphi(x_1,\ldots,x_m)$ such that $H(a_{\imag_1},\ldots,a_{\imag_m})$ if and only if $\Sa M \models \varphi(a_{\imag_1},\ldots,a_{\imag_m})$ for any $(\imag_1,\ldots,\imag_m)$ in $\omega^m$.
Finally, let $X_p$ be the set of $((d_1,e_1),\ldots,(d_m,e_m)) \in M^{2m}$ with $\Sa M \models \varphi(d_1,\ldots,d_m)$.
\end{proof}




We now consider examples.
The most basic example is the model companion of the theory of a set equipped with an $n$-ary relation for each $n \ge 1$.
Trace maximality of this theory follows from Lemma~\ref{lem:max} and it is supersimple by~\cite[Thm.~B.1]{jera}.

\begin{proposition}
\label{prop:bool}
Every infinite boolean algebra is trace maximal.
\end{proposition}

A subset $A$ of a boolean algebra is \textbf{independent} if the subalgebra generated by $A$ is free, equivalently for any $\alpha_1,\ldots,\alpha_k,\beta_1,\ldots,\beta_\ell \in A$ with $\alpha_i \ne \beta_j$ for all $i,j$ we have 
\[(\alpha_1 \land \cdots \land \alpha_k) \land ( \neg \beta_1 \land \cdots \land \neg \beta_\ell) \ne 0\]

\begin{lemma}
\label{lem:boolean}
Suppose that $\Sa B$ is a boolean algebra, $A$ is an independent subset of $B$, and $(\beta_i : i \le k)$, $(\alpha^i_j : i \in \{1,\ldots,n\}, 1 \le j \le k)$ are elements of $A$ such that for all $i \in \{1,\ldots,n\}$:
\begin{enumerate}
\item $|\{ \beta_1,\ldots,\beta_k\}| = k = |\{\alpha^i_1,\ldots,\alpha^i_k\}|$, and
\item $\{ \beta_1,\ldots,\beta_k \} \ne \{\alpha^i_1,\ldots,\alpha^i_k\}$. 
\end{enumerate}
Then $(\beta_1 \land \cdots \land \beta_k) \nleqslant \bigvee_{i = 1}^{n} (\alpha^i_1 \land \cdots \land \alpha^i_k)$.
\end{lemma}

Here and below we let $\leqslant$ be the boolean order.

\begin{proof}
Note that for each $i \in \{1,\ldots,n\}$ there is $j(i)$ such that $\alpha^i_{j(i)} \notin \{\beta_1,\ldots,\beta_k \}$.
By independence we have
$$ \gamma := (\beta_1\land \cdots \land \beta_k) \land (\neg \alpha^1_{j(1)} \land \cdots \land \neg \alpha^n_{j(n)}) \ne 0. $$
Then $\gamma \leqslant (\beta_1\land \cdots \land \beta_k)$ and $\gamma \land \bigvee_{i = 1}^{n} (\alpha^i_1 \land \cdots \land \alpha^i_k) = 0$.
\end{proof}

We now prove Proposition~\ref{prop:bool}.

\begin{proof}[Proof of Proposition~\ref{prop:bool}]
Any finite boolean algebra embeds into $\Sa B$.
Hence $\Sa B$ contains an independent subset of cardinality $n$ for every $n$.
We may suppose $\Sa B$ is $\aleph_1$-saturated.
By saturation there is a countably infinite independent subset $A$ of $B$.
By Lemma~\ref{lem:max} it suffices to suppose $E$ is a $k$-hypergraph on $A$ and produce  definable $Y \subseteq B^k$ such that we have $ E(a_1,\ldots,a_k)$ if and only if $(a_1,\ldots,a_k) \in Y$ for all $a_1,\ldots,a_k \in A$.
Let $f \colon A^k \to B$ be given by $f(a_1,\ldots,a_k) = a_1 \land \cdots \land a_k$.
Let $D$ be the set of $(a_1,\ldots,a_k) \in A^k$ such that $|\{a_1,\ldots,a_k\}| = k$.
Let $(a^i : i < \omega)$ be an enumeration of $E$.
Lemma~\ref{lem:boolean} shows that if $b \in D$ and $\neg E(b)$ then $f(b) \nleqslant f(a^1)\vee\cdots\vee f(a^n)$ for all $n$.
Hence for any $n$ there is a $c \in B$ such that $f(a^i) \le c$ for all $i \leqslant n$ and $f(b) \nleqslant c$ for all $b \in D$ such that $\neg E(b)$.
By saturation there is $c \in B$ such that $f(a^i) \leqslant c$ for all $i < \omega$ and $f(b) \nleqslant c$ for all $b \in D$ with $\neg E(b)$.
Let $Y$ be the set of $a \in B^k$ such that $a \in D$ and $f(a) \leqslant c$.
\end{proof}

I do not know of an $\ip$ field which is not trace maximal.
It is easy to see that the ring $\Z$ is trace maximal\footnote{One can use the standard codings, or alternatively apply Proposition~\ref{prop:bool}.}, and it follows that the main examples of logically wild fields, e.g. finitely generated infinite fields, are trace maximal as they interpret $\Z$.
We show below that some of the main examples of logically tame $\ip$ fields are trace maximal.
Let $K$ be a field and let $\chara(K)$ be the characteristic of $K$.
When $\chara(K) = p$ we let $\wp \colon K \to K$ be the Artin-Schreier map $\wp(x) = x^p - x$.
Recall that $K$ is {\bf Artin-Schreier closed} either $\chara(K) = 0$ or $K$ is positive characteristic and $\wp(K) = K$.

\begin{proposition}\label{prop:Fp((t))}
Suppose that $K$ is a field which is not Artin-Schreier closed.
Then $K$ is trace maximal.
\end{proposition}

It follows in particular that positive characteristic local fields are trace maximal.
Proposition~\ref{prop:Fp((t))} essentially follows from the proof of Hempel's result that a non-Artin-Schreier closed field is $k$-$\ip$ for every $k \ge 1$~\cite{hempel-field}.
Fact~\ref{fact:hempel} is the key tool used in her proof~\cite[Thm.~6.3]{hempel-field}.

\begin{fact}\label{fact:hempel}
Suppose that $K$ is not Artin-Schreier closed and let $A, B$ be disjoint finite subsets of $K$ such that $A \cup B$ is linearly independent over the prime subfield.
Then there is $\gamma \in K$ such that $a\gamma$ is in $\wp(K)$ for any $a \in A$ and $b\gamma$ is not in $\wp(K)$ for any $b \in B$.
\end{fact}

\begin{proof}[Proof of Proposition~\ref{prop:Fp((t))}]
After possibly passing to an elementary extension we suppose that $K$ is $\aleph_1$-saturated.
Fix a countably infinite subset $A$ of $K$ which is algebraically independent over the prime subfield.
By Lemma~\ref{lem:max} it is enough to suppose that $E$ is an $m$-hypergraph on $A$ and produce definable $Y \subseteq K^m$ so that we have $E(\alpha_1,\ldots,\alpha_m)$ if and only if $(\alpha_1,\ldots,\alpha_m) \in Y$ for all $\alpha_1,\ldots,\alpha_m \in A$.
Let $D$ be the set of $(\alpha_1,\ldots,\alpha_m) \in A^m$ such that $|\{ \alpha_1,\ldots,\alpha_m \}| = m$ and $f$ be the map $D \to A$ given by $f(\alpha_1,\ldots,\alpha_m) = \alpha_1 \alpha_2 \ldots \alpha_m$.
By algebraic independence we have $f(\alpha_1, \ldots, \alpha_m) = f(\beta_1, \ldots, \beta_m)$ if and only if $\{\alpha_1,\ldots,\alpha_m\} = \{\beta_1,\ldots,\beta_m\}$ for all $(\alpha_1,\ldots,\alpha_m), (\beta_1,\ldots,\beta_m) \in D$.
Furthermore $\{ f(\alpha) : \alpha \in D \}$ is linearly independent over the prime subfield.
By Fact~\ref{fact:hempel} and saturation there is $\gamma \in K$ such that we have $E(\alpha)$ if and only if $K \models \exists x\left[ \wp(x) = \gamma f(\alpha) \right]$ for all $\alpha \in D$.
\end{proof}

We now consider $\pac$ and $\prc$ fields.
See \cite{field-arithmetic, Chatzidakis} for background on $\pac$ fields.
Psuedofinite fields and infinite algebraic extensions of finite fields are $\mathrm{PAC}$ (and also simple~\cite{Chatzidakis}).

\begin{proposition}
\label{prop:pseudofinite}
Let $K$ be a field.
Suppose that $K$ is not separably closed or real closed and that some finite extension of $K$ is $\pac$.
Then $K$ is trace maximal.
\end{proposition}

This applies to $\prc$ fields that are neither separably nor real closed as $K(\sqrt{-1})$ is $\pac$ when $K$ is $\prc$.
We make use of Duret's proof~\cite{duret} that a non-separably closed $\mathrm{PAC}$ field is $\ip$.
Fact~\ref{fact:macintyre}, a Galois-theoretic exercise, was essentially proven by Macintryre~\cite{Macintyre-omegastable}.

\begin{fact}
\label{fact:macintyre}
Suppose that $K$ is not separably closed.
Then there is a finite extension $F$ of $K$ such that one of the following holds:
\begin{enumerate}
\item $F$ is not Artin-Schreier closed, or
\item there is a prime $q \ne \chara(K)$ such that $F$ contains a primitive $q$th root of unity and some element of $F$ is not a $q$th power.
\end{enumerate}
\end{fact}

Fact~\ref{fact:duret} is~\cite[Lemma~6.2]{duret}.

\begin{fact}
\label{fact:duret}
Suppose that $K$ is a $\mathrm{PAC}$ field, $q \ne \mathrm{Char}(K)$ is a prime, $K$ contains a primitive $q$th root of unity, and some element of $K$ is not a $q$th power.
Let $A,B$ be finite disjoint subsets of $K$.
Then there is $\gamma \in K$ such that $a+\gamma$ is a $q$th power for every $a \in A$ and $b+\gamma$ is not a $q$th power for every $b \in B$.
\end{fact}



\begin{proof}[Proof of Proposition~\ref{prop:pseudofinite}.]
It suffices to show that some finite extension of $K$ is trace maximal as any finite extension is interpretable in $K$.
Now $K$ has a finite $\pac$ extension which is not real or separably closed as $K$ is not real or separable closed.
Hence we may suppose that $K$ is $\pac$.
By Fact~\ref{fact:macintyre} there is a finite field extension $F$ of $K$ such that:
\begin{enumerate}
\item $F$ is not Artin-Schreier closed, or
\item there is a prime $q \ne \chara(K)$ such that $F$ contains a primitive $q$th root of unity and some element of $F$ is not a $q$th power.
\end{enumerate}
If (1) holds then $F$ is trace maximal by Proposition~\ref{prop:Fp((t))}.
Suppose $(2)$.
Now $F$ is $\mathrm{PAC}$ as $\mathrm{PAC}$ fields are closed under finite extensions.
Again, $F$ is not separably closed as $K$ is not real or separably closed.
Hence we suppose that $K$ satisfies (2) and fix a relevant prime $q$.
After possibly passing to an elementary extension we also suppose that $K$ is $\aleph_1$-saturated.
Let $A$ be a countably infinite subset of $K$ and $t$ be an element of $K$ which is not in the algebraic closure of the subfield generated by $A$.
Let $f \colon A^m \to K$ be given by $f(\alpha_0,\ldots,\alpha_{m - 1}) = t^m + \alpha_{m - 1} t^{m - 1} + \cdots + \alpha_1 t + \alpha_0$.
Let $X$ be a subset of $A^m$.
By injectivity of $f$, Fact~\ref{fact:duret}, and saturation there is $\gamma \in K$ such that we have  $\alpha \in X$ if and only if $K \models \exists x(f(\alpha) + \gamma=x^q)$ for all $\alpha \in A^m$.
By Lemma~\ref{lem:max} $K$ is trace maximal.
\end{proof}

\bibliographystyle{abbrv}
\bibliography{NIP}
\end{document}